\documentclass[11pt,reqno]{amsart}
\usepackage{enumerate}
\usepackage[margin=1in]{geometry}
\usepackage{amssymb}
\usepackage{amsmath}
\usepackage{amsthm}
\usepackage{float}
\usepackage{microtype}
\usepackage{mathtools}
\usepackage[hidelinks]{hyperref}
\usepackage{enumitem}
\usepackage{tikz-cd}

\def\quotes#1{``#1''}
\def\o#1{\operatorname{#1}}
\def\p#1{\left( #1 \right)}
\def\set#1{\left\{ #1 \right\}}
\def\abs#1{\left| #1 \right|}

\def\F{\mathbb{F}}

\def\Z{\mathbb{Z}}
\def\Q{\mathbb{Q}}
\def\R{\mathbb{R}}

\def\GL{\operatorname{GL}}
\def\SL{\operatorname{SL}}
\def\PGL{\operatorname{PGL}}
\def\PSL{\operatorname{PSL}}
\def\Cs{\mathtt{Cs}}
\def\Ns{\mathtt{Ns}}
\def\Cn{\mathtt{Cn}}
\def\Nn{\mathtt{Nn}}
\def\B{\mathtt{B}}

\def\Afour{\mathtt{A4}}
\def\Afive{\mathtt{A5}}
\def\Sfour{\mathtt{S4}}
\def\DPGL{\mathtt{GL}}
\def\DPSL{\mathtt{SL}}
\def\la#1{\mathtt{#1}}

\def\Gal{\operatorname{Gal}}
\def\Aut{\operatorname{Aut}}
\def\tr{\operatorname{tr}}
\def\Frob{\operatorname{Frob}}
\def\fp{\mathfrak{p}}
\def\leg#1#2{\left(\frac{#1}{#2}\right)}
\def\tors{\text{tors}}
\def\disc{\o{disc}}

\theoremstyle{plain}
\newtheorem{theorem}{Theorem}[section]
\newtheorem{lemma}[theorem]{Lemma}
\newtheorem{proposition}[theorem]{Proposition}

\theoremstyle{definition}

\newtheorem{remark}[theorem]{Remark}

\allowdisplaybreaks

\author[Jacob Mayle]{Jacob Mayle}
\address{Department of Mathematics,
Wake Forest University,
Winston-Salem, NC, 27104 USA }
\email{maylej@wfu.edu}
\date{}

\begin{document}

\title[Rigidity in Elliptic Curve Local-Global Principles]{Rigidity in Elliptic Curve Local-Global Principles}

\begin{abstract}
We study the rigidity of the local conditions in two well-known local-global principles for elliptic curves over number fields. In particular, we consider a local-global principle for torsion due to Serre and Katz, and one for isogenies due to Sutherland. For each of these local-global principles, we prove that if an elliptic curve $E$ over a number field $K$ is such that it fails to satisfy the local condition for at least one prime ideal of $K$ of good reduction, then $E$ can satisfy the local condition at no more than 75\% of prime ideals. We also give for (conjecturally) all elliptic curves over the rationals without complex multiplication, the densities of primes that satisfy the local conditions mentioned above.
\end{abstract}



\maketitle

\section{Introduction} \label{intro}
Let $E$ be an elliptic curve over a number field $K$. For a prime ideal $\fp \subseteq \mathcal{O}_K$ of good reduction for $E$, we write $E_\fp$ to denote the reduction of $E$ modulo $\fp$. We say that a property holds for $E$ \textit{locally everywhere} if it holds for each reduced elliptic curve $E_\fp$. Given a property that holds locally everywhere, it is natural to ask if some corresponding property holds globally. If so, such an implication is referred to as a \textit{local-global} principle. 

Two well-known local-global principles address the following questions. Fix a prime number $\ell$.
\begin{enumerate}[label=(\Alph*)]
  \item \label{K-q} If $E$ has nontrivial rational $\ell$-torsion locally everywhere, must $E$ have nontrivial rational $\ell$-torsion?
  \item \label{S-q} If $E$ admits a rational $\ell$-isogeny locally everywhere, must $E$ admit a rational $\ell$-isogeny?
\end{enumerate}
It turns out that the answer to both of these questions is \quotes{no}. For instance, consider the elliptic curves over $\Q$ with LMFDB \cite{LMFDB} labels \href{https://www.lmfdb.org/EllipticCurve/Q/11/a/1}{\texttt{11.a1}} and \href{https://www.lmfdb.org/EllipticCurve/Q/2450/i/1}{\texttt{2450.i1}}. They are given by the minimal Weierstrass equations
\begin{align*}
E_{\mathtt{11.a1}}&: y^2 + y = x^3 - x^2 - 7820x - 263580, \\ 
E_{\mathtt{2450.i1}}&: y^2 + xy = x^3 - x^2 - 107x - 379.
\end{align*}
These curves provide counterexamples for the above questions. The curve $E_{\mathtt{11.a1}}$ has nontrivial rational $5$-torsion at every prime of good reduction but has trivial rational 5-torsion itself. The curve $E_{\mathtt{2450.i1}}$ admits a rational $7$-isogeny at every prime of good reduction, but does not admit a rational $7$-isogeny itself.

Further analysis of the above counterexamples reveals some structure. One sees that $E_{\mathtt{11.a1}}$ is isogenous to an elliptic curve $E_{\mathtt{11.a1}}'$ that has nontrivial rational $5$-torsion. In addition,  the curve $E_{\mathtt{2450.i1}}$ admits a 7-isogeny not over $\Q$, but over $\Q(\sqrt{-7})$. With these examples in mind, there is the prospect of powerful local-global principles stemming from questions \ref{K-q} and \ref{S-q}.

Indeed, Serge Lang proposed and Katz proved a local-global principle that corresponds to \ref{K-q} and deals more generally with composite level. In its statement below, we write \quotes{locally almost everywhere} to mean a mildly relaxed variant of \quotes{locally everywhere} in which the corresponding local condition is only asserted to hold for a set of prime ideals of density one, in the sense natural density (which we shall define shortly). 

\begin{theorem}[Katz, 1981 \cite{Ka1981}] \label{K-thrm} Fix an integer $m \geq 2$. If the condition $\abs{E_\fp(\mathcal{O}_K/\fp)} \equiv 0\pmod m$ holds locally almost everywhere, then $E$ is $K$-isogenous to an elliptic curve $E'/K$ for which $\abs{E'_{\tors}(K)} \equiv 0 \pmod m$. 
\end{theorem}

In fact, we shall only consider the case where $m$ is prime, which dates back to two exercises of Serre \cite[p. I-2 and p. IV-6]{MR1484415}. Much more recently, Sutherland established a local-global principle associated with \ref{S-q}.
\begin{theorem}[Sutherland, 2012 \cite{Su2012}] \label{S-thrm} Fix a prime number $\ell$ for which $\sqrt{\leg{-1}{\ell}\ell} \not\in K$. Suppose that the condition that $E_\fp$ admits an $\mathcal{O}_K/\fp$-rational $\ell$-isogeny holds locally almost everywhere. Then there exists a quadratic extension $L/K$ such that then $E$ admits an $L$-rational $\ell$-isogeny. Further, if $\ell = 2, 3$ or $\ell \equiv 1 \pmod 4$, then in fact $E$ admits a $K$-rational $\ell$-isogeny.
\end{theorem}

Sutherland's result sparked an outpouring of research. Notably, Anni \cite{An2014} showed that in Theorem \ref{S-thrm}, $L$ may be taken to be $K(\sqrt{-\ell})$ and gave an explicit upper bound (depending on $K$) on the prime numbers $\ell$ for which there exists an elliptic curve $E/K$ that admits a rational $\ell$-isogeny locally everywhere, but not globally. Vogt \cite{Vo2018} gave an extension of Sutherland's result to composite level. Other authors made contributions as well, such as Banwait-Cremona \cite{BC2014} and Etropolski \cite{Et2015}. Recently, there has been increased interest in probabilistic local-global principles \cite{MR4450609} and extensions to higher dimensional abelian varieties \cite{MR4301391, MR2314732}.

In this paper, we prove that for a given elliptic curve over a number field and prime number $\ell$, a failure of either of the  \quotes{locally everywhere} conditions of \ref{K-q} or \ref{S-q} must be fairly substantial. This phenomenon is a consequence of the properties of the general linear group $\GL_2(\ell)$, as we shall see. Moreover, it contrasts the elliptic curve local-global principles of Katz and Sutherland with, for instance, the familiar local-global principle of Hasse-Minkowski. As an example, consider the equation
\begin{equation} \label{limited-ex} x^2 + y^2 = 3. \end{equation}
It fails to have solutions locally everywhere and hence has no solutions over $\Q$. However, the failure is quite limited. In fact, (\ref{limited-ex}) has no solutions over $\Q_2$ and $\Q_3$ but has solutions over $\R$ and $\Q_p$ for each $p \geq 5$.

To discuss this feature of the prime-level local-global principles of Katz and Sutherland more precisely, we fix some standard notation and terminology from algebraic number theory. Let $\mathcal{O}_K$ denote the ring of integers of $K$ and let $\mathcal{P}_K$ denote the set of prime ideals of $\mathcal{O}_K$. For a prime ideal $\fp \in \mathcal{P}_K$, denote its \textit{residue field} by $\F_\fp \coloneqq \mathcal{O}_K / \fp$ and its \textit{norm} by $N\fp \coloneqq \abs{\F_\fp}$. For a subset $\mathcal{A} \subseteq \mathcal{P}_K$ and a positive real number $x$, define
\[ \mathcal{A}(x) \coloneqq \set{\fp \in \mathcal{A}: N \fp \leq x}. \]
The \textit{natural density} of $\mathcal{A}$ is defined to be the following limit (provided it exists),
\begin{equation}
\delta(\mathcal{A}) \coloneqq \lim_{x \to \infty} \frac{\abs{\mathcal{A}(x)}}{\abs{\mathcal{P}_K(x)}}, \label{density}
\end{equation}

The subsets of $\mathcal{P}_K$ that are relevant to our study are the following,
\begin{align}
\mathcal{S}^1_{E,\ell} &\coloneqq \set{\fp \in \mathcal{P}_K : \fp \nmid N_E \text{ and } E_\fp \text{ has an } \F_\fp\text{-rational point of order } \ell}, \label{K-set} \\
\mathcal{S}_{E,\ell} &\coloneqq \set{\fp \in \mathcal{P}_K : \fp \nmid N_E \text{ and } E_\fp \text{ has an } \F_\fp\text{-rational isogeny of degree } \ell}. \label{S-set}
\end{align}
With these sets defined and the above notation in mind, we now state our main theorem.
\begin{theorem} \label{main-thrm} Let $K$ be a number field, $E/K$ be an elliptic curve, and $\ell$ be a prime number.
\begin{enumerate}
\item \label{mt-1} If there exists a prime ideal $\fp \subseteq \mathcal{O}_K$ of good reduction for $E$ such that $E_\fp$ does not have an $\mathbb{F}_\fp$-rational point of order $\ell$, then $\delta(\mathcal{S}^1_{E,\ell}) \leq \frac{3}{4}$.
\item \label{mt-2} If there exists a prime ideal $\fp \subseteq \mathcal{O}_K$ of good reduction for $E$ such that $E_\fp$ does not admit an $\mathbb{F}_\fp$-rational isogeny of degree $\ell$, then $\delta(\mathcal{S}_{E,\ell}) \leq \frac{3}{4}$.
\end{enumerate}
Moreover, if $\ell = 2$, then the quantity $\frac{3}{4}$ may be replaced with $\frac{2}{3}$ in both parts (\ref{mt-1}) and (\ref{mt-2}) above.
\end{theorem}

We prove this theorem in Section \ref{sec-thrm} by applying the Chebotarev density theorem and Propositions \ref{S-prop} and  \ref{K-prop}, which are purely group-theoretic. We complete the proofs of the two group-theoretic propositions by considering subgroups of $\GL_2(\ell)$ case-by-case in Section \ref{S-casework} and Section \ref{K-casework}, following Dickson's well-known classification. 

Our result weakens the hypotheses of Theorem \ref{K-thrm} and Theorem \ref{S-thrm}, for each reducing the density at which the local condition of its statement must hold from $1$ down to $\frac{3}{4}$ (or $\frac{2}{3}$ in the case of $\ell = 2$). Perhaps more to the point, our result may be viewed as one about the rigidity of the \quotes{locally everywhere} conditions of \ref{K-q} and \ref{S-q}. Roughly speaking, a collection is termed \textit{rigid} if its elements are determined by less information than expected. A well-known example is the subset of complex analytic functions among all complex functions. Another example, articulated by Jones \cite{Jo2017}, is the subset of power maps among all set functions $K \to K$, for a Galois number field $K$. In our case, for an odd prime $\ell$, the two parts of Theorem \ref{main-thrm} are equivalent to the assertion that (1') $E/K$ has nontrivial rational $\ell$-torsion locally everywhere if and only if $\delta(\mathcal{S}^1_{E,\ell}) > \frac{3}{4}$, and (2') $E/K$ admits a rational $\ell$-isogeny locally everywhere if and only if $\delta(\mathcal{S}_{E,\ell}) > \frac{3}{4}$. In this sense, for a number field $K$ and prime number $\ell$, the subset of elliptic curves over $K$ that satisfy the \quotes{locally everywhere} condition of \ref{K-q} (respectively \ref{S-q}) is rigid among the set of all elliptic curves over $K$.

The related matter of computing the densities $\delta(\mathcal{S}^1_{E,\ell})$ and $\delta(\mathcal{S}_{E,\ell})$ is straightforward in light of \cite{Su2016}. As we shall see, for a given elliptic curve $E/K$ and prime number $\ell$, these densities are determined  by the image of the mod $\ell$ Galois representation of $E$ in $\GL_2(\ell)$. In Section \ref{appendix}, we list the values of $\delta(\mathcal{S}^1_{E,\ell})$ and $\delta(\mathcal{S}_{E,\ell})$ corresponding to all 63 of the known (and conjecturally all) mod $\ell$ Galois images of elliptic curves over the rationals without complex multiplication.

\section{Preliminaries on Galois Representations}

In this section, we recall some basic facts about Galois representations of elliptic curves. Let $E$ be an elliptic curve over a perfect field $K$. Let $\overline{K}$ be an algebraic closure of $K$ and let $\ell$ be a prime number.

The $\ell$-torsion subgroup of $E(\overline{K})$, denoted $E[\ell]$, is a $\Z/\ell\Z$-vector space of rank two. The absolute Galois group $G_K \coloneqq \Gal(\overline{K}/K)$ acts coordinate-wise on $E[\ell]$. This action is encoded in the group homomorphism
\[ \begin{tikzcd}[sep = scriptsize]
\rho_{E,\ell}: G_K \arrow{r} & \Aut(E[\ell]) \arrow["\sim"]{r} & \GL_2(\ell),
\end{tikzcd} \]
which is known as the \textit{mod $\ell$ Galois representation} of $E$. Above $\GL_2(\ell)$ denotes the general linear group over $\F_\ell$ and the isomorphism $\Aut(E[\ell]) \overset{\sim}{\to} \GL_2(\ell)$ is determined by a choice of $\Z/\ell\Z$-basis of $E[\ell]$. The \textit{mod $\ell$ Galois image} of $E$, denoted $G_E(\ell)$, is the image of $\rho_{E,\ell}$. Because $\rho_{E,\ell}$ and $G_E(\ell)$ depend on a choice of basis for $E[\ell]$, we recognize that we may only speak sensibly of these objects up to conjugation in $\GL_2(\ell)$. 

Let $K(E[\ell])$ denote the \textit{$\ell$-division field} of $E$, that is, the Galois extension of $K$ obtained by adjoining to $K$ the affine coordinates of the points of $E[\ell]$. Observe that $\Gal(\overline{K}/K(E[\ell]))$ is the kernel of $\rho_{E,\ell}$. Thus, by the first isomorphism theorem and Galois theory,
\[ \begin{tikzcd}[sep = scriptsize] \tilde{\rho}_{E,\ell}:  \Gal(K(E[\ell])/K) \arrow["\sim"]{r} & G_E(\ell)  \end{tikzcd} \]
is an isomorphism, where $\tilde{\rho}_{E,\ell}$ is the restriction of $\rho_{E,\ell}$ to $\Gal(K(E[\ell])/K)$.

The Galois image $G_E(\ell) \subseteq \GL_2(\ell)$ is of central interest to us since it detects the presence of nontrivial rational $\ell$-torsion of $E$ and rational $\ell$-isogenies admitted by $E$. We shall describe precisely how in the following lemma. First, recall that the \textit{Borel subgroup} and \textit{first Borel subgroup} of $\GL_2(\ell)$ are, respectively,
\begin{align} \label{borel-def}
\mathcal{B}(\ell) &\coloneqq \set{\begin{pmatrix} a & b \\ 0 & d \end{pmatrix} : a,d \in \F_\ell^\times \text{ and } b \in \F_\ell},  \\
\mathcal{B}_1(\ell) &\coloneqq \set{\begin{pmatrix} 1 & b \\ 0 & d \end{pmatrix} : d \in \F_\ell^\times \text{ and } b \in \F_\ell}.
\end{align}

\begin{lemma} \label{GalImage-L} With the notation above, we have that 
\begin{enumerate}
  \item $E$ has nontrivial $K$-rational $\ell$-torsion if and only if $G_E(\ell)$ is conjugate to a subgroup of $\mathcal{B}_1(\ell)$,
  \item $E$ admits a $K$-rational $\ell$-isogeny if and only if $G_E(\ell)$ is conjugate to a subgroup of $\mathcal{B}(\ell)$.
\end{enumerate}
\end{lemma}
\begin{proof}
\begin{enumerate}
  \item Let $P \in E(K)$ be a point of order $\ell$. Then $P \in E[\ell]$ and we may choose a point $Q \in E[\ell]$ such that $\set{P,Q}$ is a $\Z/\ell\Z$-basis of $E[\ell]$. For each $\sigma \in G_K$, we have that
  \[    \sigma(P) = P \quad \text{and} \quad \sigma(Q) = b P + d Q  \]
  for some $b, d \in \Z/\ell\Z$ (depending on $\sigma$). Hence,
  \[
    \rho_{E,\ell}(\sigma) = \begin{pmatrix} 1 & b \\ 0 & d \end{pmatrix} \in \mathcal{B}_1(\ell).
  \]
  Thus $G_E(\ell) \subseteq \mathcal{B}_1(\ell)$, with respect to the basis $\set{P,Q}$. 

  Conversely, assume that $G_E(\ell)$ is conjugate to a subgroup of $\mathcal{B}_1(\ell)$. Let $\set{P,Q}$ be a $\Z/\ell\Z$-basis of $E[\ell]$ that realizes $G_E(\ell) \subseteq \mathcal{B}_1(\ell)$. Then $\sigma(P) = P$ for each $\sigma \in G_K$, so $P \in E(K)$. Thus $P$ is a nontrivial $K$-rational $\ell$-torsion point of $E$.

  \item Let $\phi: E \to E'$ be a $K$-rational $\ell$-isogeny. Then $\ker \phi \subseteq E(K)$ is cyclic of order $\ell$. Let $P$ be a generator of $\ker \phi$. Then $P \in E[\ell]$ and we may choose a point $Q \in E[\ell]$ such that $\set{P,Q}$ is a $\Z/\ell\Z$-basis of $E[\ell]$. For each $\sigma \in G_K$, we have that
  \[
    \sigma(P) = a P \quad \text{and} \quad \sigma(Q) = b P + d Q
  \]
  for some $a,b,d \in \Z/\ell\Z$ (depending on $\sigma$). Hence,
  \[
    \rho_{E,\ell}(\sigma) = \begin{pmatrix} a & b \\ 0 & d \end{pmatrix} \in \mathcal{B}(\ell).
  \]
  Thus $G_E(\ell) \subseteq \mathcal{B}(\ell)$, with respect to the basis $\set{P,Q}$. 

  Conversely, assume that $G_E(\ell)$ is conjugate to a subgroup of $\mathcal{B}(\ell)$. Let $\set{P,Q}$ be a $\Z/\ell\Z$-basis of $E[\ell]$ that realizes $G_E(\ell) \subseteq \mathcal{B}(\ell)$. Let $\Phi$ denote the subgroup of $E(K)$ generated by $P$. We have that $\Phi$ is cyclic of order $\ell$, so the natural isogeny $E \to E/\Phi$ is a $K$-rational $\ell$-isogeny of $E$. \qedhere
\end{enumerate}
\end{proof}

\section{Preliminaries on \texorpdfstring{$\GL_2(\ell)$}{GL2(l)}} \label{S-GL2}

Let $\ell$ be an odd prime number. The main objective of this section is to state two important classification results for $\GL_2(\ell)$: the classification of its subgroups (originally due to Dickson) and the classification of its conjugacy classes. We start by recalling some standard notation and terminology.

We write $\mathcal{Z}(\ell)$ to denote the center of $\GL_2(\ell)$, which consists precisely of the scalar matrices of $\GL_2(\ell)$. The \textit{projective linear group} over $\F_\ell$ is the quotient of $\PGL_2(\ell) \coloneqq \GL_2(\ell)/\mathcal{Z}(\ell)$ and $\pi: \GL_2(\ell) \twoheadrightarrow \PGL_2(\ell)$ denotes the quotient map. We denote the image of a matrix $\gamma \in \GL_2(\ell)$ in $\PGL_2(\ell)$ by $\overline{\gamma}$. Similarly, we denote the image of a subset $S \subseteq \GL_2(\ell)$ in $\PGL_2(\ell)$ by $\overline{S}$. In particular, the \textit{projective special linear group} over $\F_\ell$ is $\PSL_2(\ell) \coloneqq \overline{\SL_2(\ell)}$, where $\SL_2(\ell)$ denotes the special linear group over $\F_\ell$.

The \textit{split Cartan subgroup} of $\GL_2(\ell)$ and its normalizer are, respectively,
\begin{align*}
\mathcal{C}_s(\ell) &\coloneqq \set{\begin{pmatrix} a & 0 \\ 0 & b \end{pmatrix} : a,b \in \F_\ell^\times }, \\
\mathcal{C}^+_s(\ell) &= \mathcal{C}_s(\ell) \cup \begin{pmatrix} 0 & 1 \\ 1 & 0 \end{pmatrix} \mathcal{C}_s(\ell).
\end{align*}
Fix a non-square $\varepsilon \in \F_\ell^\times\setminus\F_\ell^{\times 2}$. The \textit{non-split Cartan} subgroup of $\GL_2(\ell)$ and its normalizer are, respectively,
\begin{align*}
\mathcal{C}_{ns}(\ell) &\coloneqq \set{\begin{pmatrix} a & \varepsilon b \\ b & a\end{pmatrix} : a,b \in \F_\ell \text{ and } (a,b) \neq (0,0)} \\
\mathcal{C}^+_{ns}(\ell) &= \mathcal{C}_{ns}(\ell) \cup \begin{pmatrix} 1 & 0 \\ 0 & -1 \end{pmatrix} \mathcal{C}_{ns}(\ell).
\end{align*}

The Borel subgroup $\mathcal{B}(\ell)$ was defined in (\ref{borel-def}). Also, let $A_n$ and $S_n$ denote the alternating group and symmetric group, respectively, on $n$ elements. With notation set, we now state Dickson's classification \cite{Di1901}.

\begin{proposition} \label{prop-Dickson} Let $\ell$ be an odd prime and $G \subseteq \GL_2(\ell)$ be a subgroup. If $\ell$ does not divide the $\abs{G}$, then
\begin{enumerate}
  \item[$\Cs$.] $G$ is conjugate to a subgroup of $\mathcal{C}_s(\ell)$;
  \item[$\Cn$.] $G$ is conjugate to a subgroup of $\mathcal{C}_{ns}(\ell)$, but not of $\mathcal{C}_s(\ell)$;
  \item[$\Ns$.] $G$ is conjugate to a subgroup of $\mathcal{C}^+_s(\ell)$, but not of $\mathcal{C}_s(\ell)$ or $\mathcal{C}_{ns}(\ell)$;
  \item[$\Nn$.] $G$ is conjugate to a subgroup of $\mathcal{C}^+_{ns}(\ell)$, but not of $\mathcal{C}^+_{s}(\ell)$ or $\mathcal{C}_{ns}(\ell)$;
  \item[$\Afour$.] $\overline{G}$ is isomorphic to $A_4$;
  \item[$\Sfour$.] $\overline{G}$ is isomorphic to $S_4$; or
  \item[$\Afive$.] $\overline{G}$ is isomorphic to $A_5$.
\end{enumerate}
If $\ell$ divides $\abs{G}$, then
\begin{enumerate}
\item[$\B$.] $G$ is conjugate to a subgroup of $\mathcal{B}(\ell)$, but not of $\mathcal{C}_s(\ell)$;
\item[$\DPSL$.] $\overline{G}$ equals $\PSL_2(\ell)$; or
\item[$\DPGL$.] $\overline{G}$ equals $\PGL_2(\ell)$.
\end{enumerate}
\end{proposition}
\begin{proof} See, for instance, \cite[Section 2]{Se1972}. \qedhere
\end{proof}
A subgroup $G \subseteq \GL_2(\ell)$ has \textit{type} $\Cs, \Cn$, $\Ns$, etc. according to its position in the classification.

Eigenvalues will play a central role in our study. For a matrix $\gamma \in \GL_2(\ell)$ we shall, in particular, be interested only in the eigenvalues of $\gamma$ that lie in $\F_\ell$. The existence of such eigenvalues may be detected by the \textit{discriminant} of $\gamma$, by which we mean the discriminant of the characteristic polynomial of $\gamma$,
\begin{equation*} \label{Del-def} 
\Delta(\gamma) \coloneqq \disc(\det(\gamma - x I)) =  (\tr \gamma )^2 - 4 \det \gamma. 
\end{equation*}
Further, define the quadratic character $\chi: \GL_2(\ell) \twoheadrightarrow \set{0,\pm 1}$ by
\begin{equation} \label{chi-def} \chi(\gamma) \coloneqq \leg{\Delta(\gamma)}{\ell}, \end{equation}
where $\leg{\cdot}{\cdot}$ denotes the Legendre symbol. Now since $\gamma$ has an eigenvalue in $\F_\ell$ if and only if its characteristic polynomial splits over $\F_\ell$, we have that
\begin{equation}
\gamma \text{ has an eigenvalue in } \F_\ell \quad \iff \quad \chi(\gamma) \neq -1. \label{chi-egeinvals}
\end{equation}

The conjugacy classes of $\GL_2(\ell)$ are well-known (see, e.g., \cite[XVIII Section 12]{La2005} or \cite[Table 1]{Su2016}). In the table below, we list the conjugacy classes of $\GL_2(\ell)$ with the associated values of $\det$, $\tr$, $\chi$, and eigenvalues.

\begin{table}[H]
\begin{tabular}{|l|c|c|c|c|c|c|}
\hline
Representative of class & No. of classes & Size of class  & $\det$ & $\tr$ & $\chi$ & Eigenvalues \\ \hline
$\begin{pmatrix} a & 0 \\ 0 & a \end{pmatrix}$ \; $0 < a < \ell$      & $\ell - 1$                    & $1$               &  $a^2$           &  $2a$ & $0$  &  $\set{a}$   \\ \hline
$\begin{pmatrix} a & 1 \\ 0 & a \end{pmatrix}$ \; $0 < a < \ell$      & $\ell - 1$                    & $(\ell+1)(\ell-1)$           &    $a^2$         & $2a$  & $0$ &  $\set{a}$      \\ \hline
$\begin{pmatrix} a & 0 \\ 0 & b \end{pmatrix}$ \; $0 < a < b < \ell$ & $\frac{1}{2}(\ell-1)(\ell-2)$ & $\ell(\ell+1)$                &   $ab$  & $a+b$  & $1$ & $\set{a,b}$    \\ \hline
$\begin{pmatrix} a & \varepsilon b \\ b & a \end{pmatrix}$ \;  $\begin{cases}0 \leq a < \ell \\  0 < b \leq \frac{\ell-1}{2} \end{cases}$ & $\frac{1}{2}\ell(\ell-1)$    & $\ell(\ell-1)$  & $a^2 - \varepsilon b^2$    &  $2a$  & $-1$ &  $\emptyset$     \\ \hline
\end{tabular} 
\caption{Conjugacy classes of $\GL_2(\ell)$}\label{conj-table} 
\end{table}

\section{Reduction to Group Theory}

In this section, we use Lemma \ref{GalImage-L} and the Chebotarev density theorem to reduce our problem from one of arithmetic geometry to one of group theory.  Let $\ell$ be a prime number and define the subsets of $\GL_2(\ell)$,
\begin{align*}
\mathcal{I}_1(\ell) &\coloneqq \set{\gamma \in \GL_2(\ell) : \text{$\gamma$ has 1 as an eigenvalue}} \\
\mathcal{I}(\ell) &\coloneqq \set{\gamma \in \GL_2(\ell) : \text{$\gamma$ has some eigenvalue in $\F_\ell$}}.
\end{align*}
 We record a quick observation that connects the sets $\mathcal{I}_1(\ell)$ and $\mathcal{I}(\ell)$ with the subgroups $\mathcal{B}_1(\ell)$ and $\mathcal{B}(\ell)$.
\begin{lemma} \label{CycGp-L} Let $G \subseteq \GL_2(\ell)$ be a cyclic subgroup and let $\gamma$ be a generator of $G$. We have
\begin{enumerate}
  \item \label{CycGp-L-1} $G$ is conjugate to a subgroup of $\mathcal{B}_1(\ell)$ if and only if $\gamma \in \mathcal{I}_1(\ell)$,
  \item $G$ is conjugate to a subgroup of $\mathcal{B}(\ell)$ if and only if $\gamma \in \mathcal{I}(\ell)$. 
\end{enumerate}
\begin{proof}
\begin{enumerate}
  \item Suppose $G$ is conjugate to a subgroup of $\mathcal{B}_1(\ell)$. Then $\gamma$ is, in particular, conjugate to some matrix in $\mathcal{B}_1(\ell)$. Thus, since eigenvalues are invariant under conjugation, $1$ is an eigenvalue of $\gamma$, and so $\gamma \in \mathcal{I}_1(\ell)$. Conversely, assume that $\gamma \in \mathcal{I}_1(\ell)$. As 1 is an eigenvalue of $\gamma$, it must be that $\gamma$ is conjugate to some matrix in $\mathcal{B}_1(\ell)$.  Hence $G$, being generated by $\gamma$, is conjugate to a subgroup of $\mathcal{B}_1(\ell)$.
  \item Follows similarly to (1). \qedhere
\end{enumerate}
\end{proof}
\end{lemma}

For a subgroup $G \subseteq \GL_2(\ell)$, we define the proportions
\begin{align*}
\mathcal{F}_1(G) \coloneqq \frac{\abs{G \cap \mathcal{I}_1(\ell)}}{\abs{G}} \quad \text{and} \quad \mathcal{F}(G) \coloneqq \frac{\abs{G \cap \mathcal{I}(\ell)}}{\abs{G}}.
\end{align*}
More verbosely, $\mathcal{F}_1(G)$ is the proportion of matrices in $G$ that have 1 as an eigenvalue and $\mathcal{F}(G)$ is the proportion of matrices in $G$ that have some eigenvalue in $\F_\ell$. In the translation of our problem to group theory, $\mathcal{F}_1(G)$ and $\mathcal{F}(G)$ become the central objects of study. We describe precisely how in the next proposition, but first we set up some preliminaries for its proof.

Let $L/K$ be a finite extension of number fields. For a prime ideal $\fp \in \mathcal{P}_K$ that is unramified in $L/K$, we write $\Frob_\fp \in \Gal(L/K)$ to denote the Frobenius element associated with $\fp$, which is defined up to conjugation. For a conjugation-stable subset $\mathcal{C} \subseteq \Gal(L/K)$, the Chebotarev density theorem states that
\[ \delta(\set{\fp \in \mathcal{P}_K : \fp \text{ is unramified in } L/K \text{ and } \Frob_\fp \in \mathcal{C}}) = \frac{\abs{\mathcal{C}}}{\abs{\Gal(L/K)}}. \]

For an elliptic curve $E$ over a number field $K$ and a prime number $\ell$, we define the set of \textit{bad} prime ideals
\[  \mathcal{D}_{E,\ell} \coloneqq \set{\fp \in \mathcal{P}_K : \fp \text{ is ramified in } K(E[\ell])/K \text{ or } \fp \mid N_E}. \]
Let $\fp \in \mathcal{P}_K \setminus \mathcal{D}_{E,\ell}$ be a \textit{good} prime and let $E_\fp$ denote the reduction of $E$ at $\fp$. As $K(E[\ell])/K$ is unramified at $\fp$, we may consider a Frobenius element $\Frob_{\fp} \in \Gal(K(E[\ell])/K)$. The Galois group $\Gal(\F_\fp(E_\fp[\ell])/\F_\fp)$ is a finite cyclic group, generated by the image of $\Frob_\fp$. Thus $G_{E_\fp}(\ell)$ is the cyclic group generated by $\tilde{\rho}_{E,\ell}(\Frob_\fp)$, up to conjugation in $\GL_2(\ell)$.

\begin{proposition} \label{bridge-prop} Let $K$ be a number field, $E/K$ be an elliptic curve, and $\ell$ be a prime number. Write $G_E(\ell)$ to denote the mod $\ell$ Galois image of $E$. Let $\mathcal{S}^1_{E,\ell}$ and $\mathcal{S}_{E,\ell}$ be as defined in (\ref{K-set}) and (\ref{S-set}). We have that
\[ \delta(\mathcal{S}^1_{E,\ell}) = \mathcal{F}_1(G_E(\ell)) \quad \text{and} \quad
\delta(\mathcal{S}_{E,\ell}) = \mathcal{F}(G_E(\ell)). \]
\end{proposition}
\begin{proof} We define two conjugate-stable subsets of $\Gal(K(E[\ell])/K)$,
\begin{align*}
\mathcal{C}^1_{E,\ell} &\coloneqq \set{\sigma \in \Gal(K(E[\ell]) / K) : \tilde{\rho}_{E,\ell}(\sigma) \in \mathcal{I}_1(\ell)} \quad \\
\mathcal{C}_{E,\ell} &\coloneqq \set{\sigma \in \Gal(K(E[\ell]) / K) : \tilde{\rho}_{E,\ell}(\sigma) \in \mathcal{I}(\ell)}.
\end{align*}
In addition, we define two subsets of $\mathcal{P}_K$,
\begin{align*}
\mathcal{T}^1_{E,\ell} &\coloneqq \set{\fp \in \mathcal{P}_K : \fp \text{ is unramified in } K(E[\ell]) / K \text{ and } \Frob_\fp \in \mathcal{C}^1_{E.\ell}}, \\
\mathcal{T}_{E,\ell} &\coloneqq \set{\fp \in \mathcal{P}_K : \fp \text{ is unramified in } K(E[\ell]) / K \text{ and } \Frob_\fp \in \mathcal{C}_{E.\ell}}.
\end{align*}
By Lemma \ref{GalImage-L} and Lemma \ref{CycGp-L}, the sets $\mathcal{S}^1_{E,\ell}$ (resp. $\mathcal{S}_{E,\ell}$) and $\mathcal{T}^1_{E,\ell}$ (resp. $\mathcal{T}_{E,\ell}$) agree up to the finite set of bad primes $\mathcal{D}_{E,\ell}$. Thus, in particular,
\begin{equation} \label{ST-eq} \delta(\mathcal{S}^1_{E,\ell}) = \delta(\mathcal{T}^1_{E,\ell}) \quad \text{and} \quad \delta(\mathcal{S}_{E,\ell}) = \delta(\mathcal{T}_{E,\ell}). \end{equation}
Now applying the Chebotarev density theorem to $\mathcal{T}^1_{E,\ell}$ and $\mathcal{T}_{E,\ell}$, we find that
\begin{align*}
\delta(\mathcal{T}^1_{E,\ell}) &= \frac{\abs{\mathcal{C}^1_{E,\ell}}}{\abs{\Gal(K(E[\ell])/K)}} = \frac{\abs{G_E(\ell) \cap \mathcal{I}_1(\ell)}}{\abs{G_E(\ell)}} = \mathcal{F}_1(G_E(\ell)), \\
\delta(\mathcal{T}_{E,\ell}) &= \frac{\abs{\mathcal{C}_{E,\ell}}}{\abs{\Gal(K(E[\ell])/K)}} = \frac{\abs{G_E(\ell) \cap \mathcal{I}(\ell)}}{\abs{G_E(\ell)}} = \mathcal{F}(G_E(\ell)). 
\end{align*}
Combining these with (\ref{ST-eq}) completes the proof.
\end{proof}

\section{Group Theoretic Propositions and Proof of Main Theorem} \label{sec-thrm}

Proposition \ref{bridge-prop} offers us a bridge between the realms of arithmetic geometry and group theory. Given it, our main objects of study are now $\mathcal{F}_1(G)$ and $\mathcal{F}(G)$. Explicitly, our goal is to show that as $G$ varies among all subgroups of $\GL_2(\ell)$, these proportions never take on a value in the open interval $\p{\frac{3}{4},1}$ when $\ell$ is an odd prime and in $\p{\frac{2}{3},1}$ when $\ell = 2$. We start with $\ell = 2$, simply proceeding \quotes{by hand} in this case.

\begin{remark} \label{K-rmrk} By inspection of each of the six matrices of $\GL_2(2)$, we find that
\[ \mathcal{I}_1(2) = \mathcal{I}(2) = \set{
\begin{pmatrix} 1 & 0 \\ 0 & 1 \end{pmatrix},
\begin{pmatrix} 0 & 1 \\ 1 & 0 \end{pmatrix},
\begin{pmatrix} 1 & 1 \\ 0 & 1 \end{pmatrix},
\begin{pmatrix} 1 & 0 \\ 1 & 1 \end{pmatrix}}.
\]
Given this, we now compute $\mathcal{F}_1(G)$ and $\mathcal{F}(G)$ for each of the six subgroups of $\GL_2(2)$, recording our results in the table below. \vspace{-12pt}
\begin{table}[H] \renewcommand{\arraystretch}{1.5}
\begin{tabular}{|l|c|c|}  \hline
Subgroup $G$ & $\mathcal{F}_1(G)$ & $\mathcal{F}(G)$ \\ \hline
$\set{\begin{psmallmatrix}1 & 0 \\ 0 & 1 \end{psmallmatrix}}$ & $1$ & $1$ \\ \hline
$\set{\begin{psmallmatrix}1 & 0 \\ 0 & 1 \end{psmallmatrix}, \begin{psmallmatrix}0 & 1 \\ 1 & 0 \end{psmallmatrix}}$ & $1$ & $1$ \\ \hline
$\set{\begin{psmallmatrix}1 & 0 \\ 0 & 1\end{psmallmatrix}, \begin{psmallmatrix}1 & 0 \\ 1 & 1 \end{psmallmatrix}}$ & $1$ & $1$ \\ \hline
$\set{\begin{psmallmatrix}1 & 0 \\ 0 & 1\end{psmallmatrix}, \begin{psmallmatrix}1 & 1 \\ 0 & 1 \end{psmallmatrix}}$ & $1$ & $1$ \\ \hline
$\set{\begin{psmallmatrix}1 & 0 \\ 0 & 1\end{psmallmatrix}, \begin{psmallmatrix}1 & 1 \\ 1 & 0 \end{psmallmatrix}, \begin{psmallmatrix}0 & 1 \\ 1 & 1 \end{psmallmatrix}}$ & $\frac{1}{3}$ & $\frac{1}{3}$ \\\hline
$\GL_2(2)$ & $\frac{2}{3}$ & $\frac{2}{3}$ \\ \hline
\end{tabular} \vspace{-10pt}
\end{table}
\noindent From the table, we observe that if $\mathcal{F}_1(G) \neq 1$ (resp. $\mathcal{F}(G) \neq 1$), then $\mathcal{F}_1(G) \leq \frac{2}{3}$ (resp. $\mathcal{F}(G) \leq \frac{2}{3}$).
\end{remark}

For the remainder of the article, we shall focus our attention exclusively on primes $\ell > 2$. We now state our main group-theoretic propositions, which we shall prove in Sections \ref{S-casework} and \ref{K-casework}. The first proposition concerns $\mathcal{F}(G)$.

\begin{proposition} \label{S-prop}  Let $\ell$ be an odd prime and $G \subseteq \GL_2(\ell)$ a subgroup. We have that
\[
\mathcal{F}(G) = 
\begin{cases}
1 & G \text{ is of type } \Cs \text{ or } \B \\
\frac{1}{\abs{\overline{G}}} & G \text{ is of type } \Cn \\
\frac{\ell+3}{2(\ell+1)} & G \text{ is of type } \DPSL \\
\frac{\ell+2}{2(\ell+1)} & G \text{ is of type } \DPGL \\
\end{cases}
\]
and
\[
\mathcal{F}(G) \in 
\begin{cases}
\set{\frac{1}{2}, \frac{3}{4}, 1} & G \text{ is of type } \Ns \\
\set{\frac{1}{\abs{\overline{G}}}, \frac{1}{4}+ \frac{1}{\abs{\overline{G}}}, \frac{1}{2} + \frac{1}{\abs{\overline{G}}}} & G \text{ is of type } \Nn \\
\set{\frac{1}{12}, \frac{1}{3}, \frac{3}{4}, 1} & G \text{ is of type } \Afour \\
\set{\frac{1}{24}, \frac{7}{24}, \frac{3}{8}, \frac{5}{12}, \frac{5}{8}, \frac{2}{3}, \frac{3}{4}, 1} & G \text{ is of type } \Sfour \\
\set{\frac{1}{60}, \frac{4}{15}, \frac{7}{20}, \frac{5}{12}, \frac{3}{5}, \frac{2}{3}, \frac{3}{4}, 1} & G \text{ is of type } \Afive. \\
\end{cases}
\]
\noindent Further, if $G$ is of type $\Cn, \Nn, \DPSL,$ or $\DPGL$, then $\mathcal{F}(G) \leq \frac{3}{4}$. In all cases, if $\mathcal{F}(G) \neq 1$, then $\mathcal{F}(G) \leq \frac{3}{4}$.
\end{proposition} 
\begin{proof} This result is the collection of Lemmas \ref{S-BCs}, \ref{S-Cns}, \ref{S-Ns}, \ref{S-Nns}, \ref{S-PSL}, \ref{S-PGL}, \ref{S-A4}, \ref{S-S4}, \ref{S-A5}, and Remark \ref{part-rmk}.
\end{proof}

Next is our group-theoretic proposition about $\mathcal{F}_1(G)$. A quick observation reduces the number of cases that we must consider. Because $\mathcal{I}_1(\ell) \subseteq \mathcal{I}(\ell)$, we have that $\mathcal{F}_1(G) \leq \mathcal{F}(G)$. Since $\mathcal{F}(G) \leq \frac{3}{4}$ holds by the above proposition when $G$ is of type $\Cn, \Nn, \DPSL,$ or  $\DPGL$, we already have that $\mathcal{F}_1(G) \leq \frac{3}{4}$ for each of these types. Thus in the following proposition, we need only consider subgroups of type $\Cs, \Ns,\B, \Afour, \Sfour,$ and $\Afive$.

\begin{proposition} \label{K-prop} Let $\ell$ be an odd prime and $G \subseteq \GL_2(\ell)$ be a subgroup. If $\mathcal{F}_1(G) \neq 1$, then
\[
\mathcal{F}_1(G) \leq
\begin{cases}
\frac{1}{2} + \frac{1}{\abs{G}} & G \text{ is of type } \Cs \text{ or } \Ns \\
\frac{1}{2} + \frac{\ell}{\abs{G}} & G \text{ is of type } \B \\
\frac{3}{4} & G \text{ is of type } \Afour, \Sfour, \text{ or } \Afive \\
\end{cases}
\]
In all cases, if $\mathcal{F}_1(G) \neq 1$, then $\mathcal{F}_1(G) \leq \frac{3}{4}$. 
\end{proposition}
\begin{proof} This result is the collection of Lemmas \ref{proj-lem}, \ref{K-Cs}, \ref{K-B}, \ref{K-Ns}, \ref{K-exp-lem}, \ref{K-exp}, and Remark \ref{K-B-ext}.
\end{proof}

The work carried out in Section \ref{S-casework} and Section \ref{K-casework} complete the proofs of the above propositions. We now prove our main theorem, given the two propositions.

\begin{proof}[Proof of Theorem \ref{main-thrm}] We prove part (1) of the statement of the theorem, noting that (2) follows in the same way. Suppose that the condition that $E$ has nontrivial rational $\ell$-torsion locally everywhere fails. Let $\fp \in \mathcal{P}_K$ be a prime ideal of good reduction for $E$ with the property that the reduction $E_{\fp}$ has trivial $\F_\fp$-rational $\ell$-torsion. Then, by Lemma \ref{GalImage-L}, the group $G_{E_{\fp}}(\ell)$ is not conjugate to a subgroup of $\mathcal{B}_1(\ell)$. Thus by Lemma \ref{CycGp-L}(\ref{CycGp-L-1}), we have that $\tilde{\rho}_{E,\ell}(\Frob_\fp) \not\in \mathcal{I}_1(\ell)$. As a result, $G_E(\ell) \cap \mathcal{I}_1(\ell)$ is a proper subset of $G_E(\ell)$, and so $\mathcal{F}_1(G_E(\ell)) \neq 1$. Thus if $\ell$ is an odd prime, then by Propositions \ref{bridge-prop} and \ref{K-prop}, we have $\delta(\mathcal{S}^1_{E,\ell}) = \mathcal{F}_1(G_E(\ell)) \leq \frac{3}{4}$. If $\ell = 2$, then Remark \ref{K-rmrk} gives $\mathcal{F}_1(G) \leq \frac{2}{3}$ and hence $\delta(\mathcal{S}^1_{E,\ell})  \leq \frac{2}{3}$.
\end{proof}

\section{Proof of Proposition \ref{S-prop}} \label{S-casework}

In this section, we prove the lemmas that are referenced in our proof of Proposition \ref{S-prop}. We begin with several observations that will be useful at times. From here on, $\ell$ denotes an odd prime number.

\begin{lemma} \label{I-lemm} Each of the following statements holds:
\begin{enumerate}
  \item For $\gamma \in \GL_2(\ell)$, we have $\gamma \in \mathcal{I}(\ell)$ if and only if $\chi(\gamma) \neq -1$, where $\chi$ is defined in (\ref{chi-def}).
  \item For $\gamma \in \GL_2(\ell)$, we have $\gamma \in \mathcal{I}(\ell)$ if and only if $\gamma^2 \in \mathcal{I}(\ell) \setminus \mathcal{Z}_{nr}(\ell)$, where $\mathcal{Z}_{nr}(\ell) \coloneqq \set{\begin{psmallmatrix} a & 0 \\ 0 & a \end{psmallmatrix} : a \in \F_\ell^\times \setminus \F_\ell^{\times 2}}$.
  \item For $\gamma_1,\gamma_2 \in \GL_2(\ell)$, if $\overline{\gamma}_1$ is conjugate to $\overline{\gamma}_2$ in $\PGL_2(\ell)$, then $\gamma_1 \in \mathcal{I}(\ell)$ if and only if $\gamma_2 \in \mathcal{I}(\ell)$.
  \item For a subgroup $G \subseteq \GL_2(\ell)$, we have 
  \[ \mathcal{F}(G) = \frac{\abs{\overline{G \cap \mathcal{I}(\ell)}}}{\abs{\overline{G}}}.\] 
  \item For subgroups $G_1,G_2 \subseteq \GL_2(\ell)$, if $\overline{G}_1$ is conjugate to $\overline{G}_2$ in $\PGL_2(\ell)$, then $\mathcal{F}(G_1) = \mathcal{F}(G_2)$. In particular, if $G_1$ and $G_2$ are conjugate in $\GL_2(\ell)$, then $\mathcal{F}(G_1) = \mathcal{F}(G_2)$.
\end{enumerate}
\end{lemma}
\begin{proof}
\begin{enumerate}
  \item We have already seen this in (\ref{chi-egeinvals}) of Section \ref{S-GL2}.
  \item Suppose that $\gamma \in \mathcal{I}(\ell)$ and say $\lambda \in \F_\ell$ is an eigenvalue of $\gamma$. Then $\lambda^2 \in \F_\ell$ is an eigenvalue of $\gamma^2$, so $\gamma^2 \in \mathcal{I}(\ell)$. Now note that if $\gamma^2 \in \mathcal{Z}(\ell)$, then as $\lambda^2$ is an eigenvalue, we must have that $\gamma^2 = \begin{psmallmatrix} \lambda^2 & 0 \\ 0 & \lambda^2 \end{psmallmatrix} \not\in \mathcal{Z}_{nr}(\ell)$.
  We now prove the converse via its contrapositive. Suppose that $\gamma \not\in \mathcal{I}(\ell)$. From the classification of conjugacy classes of $\GL_2(\ell)$ given in Table \ref{conj-table}, we see that $\gamma$ is conjugate in $\GL_2(\ell)$ to a matrix of the form $\begin{psmallmatrix} a & b \varepsilon \\ b & a \end{psmallmatrix}$ for some  $a \in \F_\ell$ and $b \in \F_\ell^\times$. Thus $\gamma^2$ is conjugate to $\begin{psmallmatrix} a^2 + b^2 \varepsilon & 2ab \varepsilon \\ 2ab & a^2 + b^2 \varepsilon \end{psmallmatrix}$ and we may calculate
  \[ \chi(\gamma^2) = \leg{\p{2(a^2+b^2 \varepsilon)}^2 - 4(a^2-b^2 \varepsilon)^2}{\ell} = \leg{16 a^2 b^2 \varepsilon}{\ell} = -\leg{a}{\ell}^2. \]
  If $a \neq 0$, then $\chi(\gamma^2) = -1$ so $\gamma^2 \not\in \mathcal{I}(\ell)$ and we are done. On the other hand, if $a = 0$, then $\gamma^2$ is conjugate, and hence equal, to the scalar matrix $\begin{psmallmatrix} b^2 \varepsilon & 0 \\ 0 & b^2 \varepsilon \end{psmallmatrix} \in {Z}_{nr}(\ell)$ and we are done as well.
  \item Since $\overline{\gamma}_1$ is conjugate to $\overline{\gamma}_2$ in $\PGL_2(\ell)$, we have that $\gamma_0 \gamma_1 \gamma_0^{-1} = \alpha \gamma_2$ holds for some $\alpha \in \F_\ell^\times$ and $\gamma_0 \in \GL_2(\ell)$. In particular, then $\gamma_1$ has an eigenvalue in $\F_\ell$ if and only if $\gamma_2$ has an eigenvalue in $\F_\ell$.
  \item In this part, we abuse notation to let $\pi: G \twoheadrightarrow \overline{G}$ denote the restriction of $\GL_2(\ell) \twoheadrightarrow \PGL_2(\ell)$ to $G$. It follows from part (3) that for a matrix $\gamma$, either $\pi^{-1}(\overline{\gamma}) \subseteq \mathcal{I}(\ell)$ or $\pi^{-1}(\overline{\gamma}) \cap \mathcal{I}(\ell) = \emptyset$ according to whether $\gamma \in \mathcal{I}(\ell)$ or not. In addition, $\abs{\pi^{-1}(\overline{\gamma})} = \abs{\ker \pi} = \frac{\abs{G}}{\abs{\overline{G}}}$. With these observations, we calculate
  \[
    \mathcal{F}(G) 
    = \frac{1}{\abs{G}} \sum_{\overline{\gamma} \in \overline{G \cap \mathcal{I}(\ell)}}  \abs{\pi^{-1}(\overline{\gamma})} 
    = \frac{1}{\abs{\overline{G}}} \sum_{\overline{\gamma} \in \overline{G \cap \mathcal{I}(\ell)}} 1
    = \frac{\abs{\overline{G \cap \mathcal{I}(\ell)}}}{\abs{\overline{G}}}
  \]
  \item This follows from parts (3) and (4). \qedhere
\end{enumerate}
\end{proof}

The remainder of this section is devoted to proving the lemmas referenced in our proof of Proposition \ref{S-prop}. We proceed case-by-case along Dickson's classification of subgroups of $\GL_2(\ell)$. Throughout, $G$ denotes a subgroup of $\GL_2(\ell)$. By Lemma \ref{I-lemm}(5), the value of  $\mathcal{F}(G)$ is invariant on conjugating $G$ in $\GL_2(\ell)$. Thus if $G$ is of type $\Cs, \Cn, \Ns, \Nn,$ or $\B$, it suffices to assume that $G$ itself is contained in 
$\mathcal{C}_s(\ell)$, 
$\mathcal{C}_{ns}(\ell)$, 
$\mathcal{C}^+_s(\ell)$ but not $\mathcal{C}_s(\ell)$,
$\mathcal{C}^+_{ns}(\ell)$ but not $\mathcal{C}_{ns}(\ell)$,
or $\mathcal{B}(\ell)$, respectively. If $G$ is of type $\DPSL$ or $\DPGL$, it suffices to assume that $G$ is equal to $\SL_2(\ell)$ or $\GL_2(\ell)$, respectively.
When $G$ is of one of the types mentioned in this paragraph, we shall make the appropriate assumption listed here without any further mention.

\subsection{Cartan and Borel subgroups} \label{S-CsB-P}
\begin{lemma} \label{S-BCs} If $G$ is of type $\Cs$ or $\B$, then $\mathcal{F}(G) = 1$.
\end{lemma}
\begin{proof} This follows immediately since a matrix of the form $\begin{psmallmatrix} a & b \\ 0 & d \end{psmallmatrix}$ has eigenvalues $a,d \in \F_\ell$.
\end{proof}

\begin{lemma} \label{S-Cns} If $G$ is of type $\Cn$, then $\mathcal{F}(G) = \frac{1}{\abs{\overline{G}}}$.
\end{lemma}
\begin{proof} For a matrix $\gamma \coloneqq \begin{psmallmatrix} a & b\varepsilon \\ b & a \end{psmallmatrix} \in \mathcal{C}_{ns}(\ell)$, we calculate that 
\[ \chi(\gamma) = \leg{(2a)^2 - 4(a^2 - b^2 \varepsilon)}{\ell} = \leg{4b^2 \varepsilon}{\ell} = -\leg{b}{\ell}^2. \]
Thus $\gamma \in \mathcal{I}(\ell)$ if and only if $b = 0$. Hence $\mathcal{C}_{ns}(\ell) \cap \mathcal{I}(\ell) = \mathcal{Z}(\ell)$, and so
\begin{equation} \label{S-Cns-eq}
G \cap \mathcal{I}(\ell) = G \cap \mathcal{Z}(\ell).
\end{equation}
As such, $\overline{G \cap \mathcal{I}(\ell)} \subseteq \overline{\mathcal{Z}(\ell)} = \set{\overline{I}}$ so Lemma \ref{I-lemm}(4) gives that $\mathcal{F}(G) = \frac{1}{\abs{\overline{G}}}$.
\end{proof}

\subsection{Normalizers of Cartan subgroups} Here we first prove a straightforward auxiliary lemma.

\begin{lemma} \label{D-lemm} If $\gamma \in \mathcal{C}_s(\ell)$ (resp. $\gamma \in \mathcal{C}_{ns}(\ell)$) and $\gamma_0 \in \mathcal{C}^+_s(\ell) \setminus \mathcal{C}_s(\ell)$ (resp. $\gamma_0 \in \mathcal{C}^+_{ns}(\ell) \setminus \mathcal{C}_{ns}(\ell)$), then
\[ \Delta(\gamma \gamma_0) = \det(\gamma) \Delta(\gamma_0) . \]
\end{lemma}
\begin{proof} For given matrices $\gamma \coloneqq \begin{psmallmatrix} a & 0 \\ 0 & d \end{psmallmatrix} \in \mathcal{C}_s(\ell)$ and $\gamma_0 \coloneqq \begin{psmallmatrix} 0 & b \\ c & 0 \end{psmallmatrix}  \in \mathcal{C}^+_s(\ell) \setminus \mathcal{C}_s(\ell)$, we calculate that
\[ \Delta(\gamma \gamma_0)  = 4abcd  = \det(\gamma) \Delta(\gamma_0) . \]
Second, for given matrices $\gamma \coloneqq \begin{psmallmatrix} a & b \varepsilon \\ b & a \end{psmallmatrix}\in \mathcal{C}_{ns}(\ell)$ and $\gamma_0 \coloneqq  \begin{psmallmatrix} c & d \varepsilon \\ -d & -c \end{psmallmatrix} \in \mathcal{C}^+_{ns}(\ell) \setminus \mathcal{C}_{ns}(\ell)$, we calculate that
\[ \Delta(\gamma \gamma_0) =4(a^2 - b^2 \varepsilon)(c^2 - d^2 \varepsilon) = \det(\gamma) \Delta(\gamma_0), \]
concluding the proof.
\end{proof}

\begin{lemma} \label{S-Ns} If $G$ is of type $\Ns$, then $\mathcal{F}(G) \in \set{\frac{1}{2}, \frac{3}{4},1}$.
\end{lemma}
\begin{proof} Write $G_c \coloneqq G \cap \mathcal{C}_s(\ell)$ and $G_n \coloneqq G \setminus G_c$. We are assuming that $G_n \neq \emptyset$, so we may fix a matrix $\gamma_0 \in G_n$. Right multiplication by $\gamma_0$ gives a bijection $G_c \to G_n$. Thus $\abs{G_c}=\abs{G_n}=\frac{1}{2}\abs{G}$ and in fact
\[ G_n = \set{\gamma \gamma_0 : \gamma \in G_c}. \]
Hence, by Lemmas \ref{I-lemm}(1) and \ref{D-lemm}, we have that
\begin{align*}
G_n \cap \mathcal{I}(\ell) &= \set{\gamma \gamma_0 : \gamma \in G_c \text{ and } \chi(\gamma\gamma_0) \neq -1 } \\
&= \set{\gamma \gamma_0 :  \gamma \in G_c \text{ and } \leg{\det \gamma}{\ell} = \chi(\gamma_0) }.
\end{align*}
Noting that $\chi(\gamma_0)$ is fixed and $\leg{\det \cdot}{\ell} : G_c \to \set{\pm 1}$ is a homomorphism, we find that
\begin{align*}
\abs{G_n \cap \mathcal{I}(\ell)} &= \abs{\set{\gamma \in G_c : \leg{\det \gamma}{\ell} = \chi(\gamma_0)}} \\
& \in \set{0, \frac{1}{2}\abs{G_c}, \abs{G_c}}= \set{0, \frac{1}{4}\abs{G}, \frac{1}{2}\abs{G}}.
\end{align*}
So, we have that
\[
\mathcal{F}(G) = \frac{\abs{G_c \cap \mathcal{I}(\ell)} + \abs{G_n \cap \mathcal{I}(\ell)}}{\abs{G}} = \frac{\frac{1}{2}\abs{G} + \abs{G_n \cap \mathcal{I}(\ell)}}{\abs{G}} \in \set{\frac{1}{2}, \frac{3}{4}, 1},
\]
which concludes the proof.
\end{proof}

\begin{lemma} \label{S-Nns} If $G$ is of type $\Nn$, then  $\mathcal{F}(G) \in \set{ \frac{1}{\abs{\overline{G}}},  \frac{1}{4} + \frac{1}{\abs{\overline{G}}},  \frac{1}{2} + \frac{1}{\abs{\overline{G}}}}$.
\end{lemma}
\begin{proof} Write $G_c \coloneqq G \cap \mathcal{C}_{ns}(\ell)$ and $G_n \coloneqq G \setminus G_{c}$. We note that by (\ref{S-Cns-eq}), 
\begin{equation} \label{S-Nns-1}
\abs{G_c \cap \mathcal{I}(\ell)} = \abs{G \cap \mathcal{Z}(\ell)}  = \frac{\abs{{G}}}{\abs{\overline{G}}}.
\end{equation}
Since $G_n \neq \emptyset$,  we may fix a matrix $\gamma_0 \in G_n$. Right multiplication by $\gamma_0$ gives a bijection $G_c \to G_n$. Thus $\abs{G_c}=\abs{G_n}=\frac{1}{2}\abs{G}$ and $G_n = \set{\gamma \gamma_0 : \gamma \in G_c}$. As in the preceding proof, Lemmas \ref{I-lemm}(1) and \ref{D-lemm} give that
\[ G_n \cap \mathcal{I}(\ell) = \set{\gamma \gamma_0 :  \gamma \in G_c \text{ and } \leg{\det \gamma}{\ell} = \chi(\gamma_0) }.
\]
Noting that $\chi(\gamma_0)$ is fixed and $\leg{\det \cdot}{\ell} : G_c \to \set{\pm 1}$ is a homomorphism, we find that
\begin{align} \label{S-Nns-2}
\abs{G_n \cap \mathcal{I}(\ell)} &= \abs{\set{\gamma \in G_c : \leg{\det \gamma}{\ell} = \chi(\gamma_0)}} \\
&\in \set{0, \frac{1}{2}\abs{G_c}, \abs{G_c}} = \set{0, \frac{1}{4}\abs{G}, \frac{1}{2}\abs{G}}.
\end{align}
Combining (\ref{S-Nns-1}) and (\ref{S-Nns-2}), we obtain
\begin{align*}
\mathcal{F}(G) &= \frac{\abs{G_c \cap \mathcal{I}(\ell)} + \abs{G_n \cap \mathcal{I}(\ell)}}{\abs{G}} \\
&= \frac{\frac{\abs{G}}{\abs{\overline{G}}} + \abs{G_n \cap \mathcal{I}(\ell)}}{\abs{G}} \\
&\in \set{ \frac{1}{\abs{\overline{G}}},  \frac{1}{4} + \frac{1}{\abs{\overline{G}}},  \frac{1}{2} + \frac{1}{\abs{\overline{G}}}},
\end{align*}
which concludes the proof.
\end{proof}

\subsection{Subgroups containing the special linear group}
\begin{lemma} \label{S-PSL} If $G$ is of type $\DPSL$, then $\mathcal{F}(G) = \frac{\ell+3}{2(\ell+1)}$.
\end{lemma}
\begin{proof} We proceed by counting the complement of $\mathcal{I}(\ell)$ in $\SL_2(\ell)$. Referencing Table \ref{conj-table}, we see that 
\[ \SL_2(\ell) \setminus \mathcal{I}(\ell) = \bigcup_{\substack{ 0 \leq a < \ell \\ 0 < b \leq \frac{\ell-1}{2}  \\ a^2 - \epsilon b^2 \equiv 1 \!\!\! \pmod \ell}} \!\!\!\!\!\!\!\! \left[ \begin{pmatrix} a & b \varepsilon \\ b & a \end{pmatrix} \right], \]
where $[\gamma]$ denotes the $\GL_2(\ell)$-conjugacy class of $\gamma$. It is well-known (e.g. \cite[Problem 22 of Section 3.2]{Ni1991}) that
\[ \abs{\set{(a,b)\in \F_\ell \oplus \F_\ell: a^2 - \varepsilon b^2 = 1}} = \ell - \leg{\varepsilon}{\ell} = \ell + 1. \] 
Realizing that solutions to $a^2 - \varepsilon b^2 = 1$ come in pairs $(\pm x, y)$ and disregarding the pair $(\pm 1, 0)$, we obtain that $\abs{U} = \frac{1}{2}(\ell-1)$ where
\[ U \coloneqq \set{(a,b) : 0 \leq a < \ell, \, 0 < b \leq \tfrac{1}{2}(\ell-1), \text{ and } a^2 - \varepsilon b^2 \equiv 1 \pmod \ell}. \] 
Thus $\SL_2(\ell) \setminus \mathcal{I}(\ell)$ is the union of $\frac{1}{2}(\ell-1)$ conjugacy classes, each of size $\ell(\ell-1)$. Hence,
\[
\mathcal{F}(\SL_2(\ell)) = 1 - \frac{\abs{\SL_2(\ell)\setminus\mathcal{I}(\ell)}}{\abs{\SL_2(\ell)}} = 1 - \frac{\frac{1}{2}( \ell - 1) \cdot \ell(\ell-1)}{\ell(\ell^2-1)}  = \frac{\ell+3}{2(\ell+1)},
\]
which completes the proof.
\end{proof}

\begin{lemma} \label{S-PGL} If $G \subseteq \GL_2(\ell)$ is a subgroup of type $\DPGL$, then $\mathcal{F}(G) = \frac{\ell+2}{2(\ell+1)}$.
\end{lemma}
\begin{proof} We proceed by counting the complement of $\mathcal{I}(\ell)$ in $\GL_2(\ell)$. Referencing Table \ref{conj-table}, we see that 
\[ \GL_2(\ell) \setminus \mathcal{I}(\ell) = \bigcup_{\substack{ 0 \leq a < \ell \\ 0 < b \leq \frac{\ell-1}{2}  }} \left[ \begin{pmatrix} a & b \varepsilon \\ b & a \end{pmatrix} \right]. \]
Thus $\GL_2(\ell) \setminus \mathcal{I}(\ell)$ is the union of $\frac{1}{2}\ell(\ell-1)$ conjugacy classes, each of size $\ell(\ell-1)$. Hence,
\[
\mathcal{F}(\GL_2(\ell)) = 1 - \frac{\abs{\GL_2(\ell)\setminus\mathcal{I}(\ell)}}{\abs{\GL_2(\ell)}} = 1 - \frac{\frac{1}{2} \ell(\ell-1) \cdot \ell(\ell-1)}{(\ell^2-1)(\ell^2-\ell)} = \frac{\ell+2}{2(\ell+1)},
\]
which completes the proof.
\end{proof}

\begin{remark} \label{part-rmk} If $G$ is of type $\Cn, \Nn, \DPSL,$ or $\DPGL$, then $\mathcal{F}(G) \leq \frac{3}{4}$. Indeed, if $G$ is of type $\Cn$ or $\Nn$, we note that $\abs{\overline{G}} \geq 2$ or $\abs{\overline{G}} \geq 4$, respectively. The inequality  now follows directly from Lemmas \ref{S-Cns} and \ref{S-Nns}, respectively. For $G$ of type $\DPSL$ or $\DPGL$, apply Lemmas \ref{S-PSL} and \ref{S-PGL} and note that $\frac{\ell+2}{2(\ell+1)} \leq \frac{\ell+3}{2(\ell+1)} \leq \frac{3}{4}$ for all $\ell \geq 3$.
\end{remark}

\subsection{Exceptional subgroups} We first introduce some notation is useful in dealing with the exceptional subgroups, i.e., those of type $\Afour$, $\Sfour$, and $\Afive$. Let $G \subseteq \GL_2(\ell)$ be a subgroup and let $H$ be a group that is isomorphic to $\overline{G} \subseteq \PGL_2(\ell)$. Let $\phi: \overline{G} \overset{\sim}{\to} H$ be an isomorphism. For each $h \in H$, we define
\[
\delta_{h} \coloneqq 
\begin{cases}
1 & \phi^{-1}(h) \in \overline{G \cap \mathcal{I}(\ell)} \\
0 & \phi^{-1}(h) \not\in \overline{G \cap \mathcal{I}(\ell)}
\end{cases}.
\]
Let $h_1,\ldots,h_n \in H$ be representatives of the conjugacy classes of $H$. For each $i \in \set{1,\ldots,n}$, we write $[h_i]$ to denote the conjugacy class of $h_i$ in $H$. By parts (3) and (4) of Lemma \ref{I-lemm}, we  have that
\begin{equation} \label{S-ex-sum}
\mathcal{F}(G) =  \frac{\sum_i \delta_{h_i} \abs{\mathcal{C}_{h_i}}}{\abs{\overline{G}}}.
\end{equation}
We record that if $1_H$ denotes the identity element of $H$, then since $\phi^{-1}(1_H) = \overline{I} \in \overline{G \cap \mathcal{I}(\ell)}$, we have $\delta_{1_H} = 1$.

\begin{lemma} \label{S-A4} If $G \subseteq \GL_2(\ell)$ is a subgroup of type $\Afour$, then $\mathcal{F}(G) \in \set{\frac{1}{12}, \frac{1}{3}, \frac{3}{4}, 1}$.
\end{lemma}
\begin{proof} Let $\phi: \overline{G} \to A_4$ be an isomorphism and define $\delta_h$ for each $h \in A_4$ as above. The conjugacy classes of $A_4$ are $[()], [(12)(34)], [(123)],$ and $[(124)]$, of sizes 1,3,4 and 4, respectively. Observe that $(123)^2 = (132) \in [(124)]$. Hence, by Lemma \ref{I-lemm}(2,3), we have that $\delta_{(124)} = \delta_{(123)}$. Putting this information together with (\ref{S-ex-sum}),
\begin{align*} \mathcal{F}(G)  &= \frac{1 \cdot \delta_{()} + 3 \cdot \delta_{(12)(34)} + 4 \cdot \delta_{(123)} + 4 \cdot \delta_{(124)}}{12} \\
&= \frac{1 + 3 \cdot \delta_{(12)(34)} + 8 \cdot \delta_{(123)}}{12}.
\end{align*}
Iterating over all $\delta_{(12)(34)}, \delta_{(123)} \in \set{0,1}$, we obtain the desired result.
\end{proof}

\begin{lemma} \label{S-S4} If $G \subseteq \GL_2(\ell)$ is a subgroup of type $\Sfour$, then $\mathcal{F}(G) \in \set{\frac{1}{24}, \frac{7}{24}, \frac{3}{8}, \frac{5}{12}, \frac{5}{8}, \frac{2}{3}, \frac{3}{4}, 1}$.
\end{lemma}
\begin{proof} Let $\phi: \overline{G} \to S_4$ be an isomorphism and define $\delta_h$ for each $h \in S_4$ as above. The conjugacy classes of $S_4$ are $[()], [(12)], [(12)(34)], [(123)]$, and $[(1234)]$, of sizes 1,6,3,8, and 6, respectively. Observe that $(1234)^2 = (13)(24) \in [(12)(34)]$. Hence, by Lemma \ref{I-lemm}(2,3), we have that $\delta_{(12)(34)} = \delta_{(1234)}$. Thus by (\ref{S-ex-sum}),
\begin{align*}
\mathcal{F}(G) &= \frac{1 \cdot \delta_{()} + 6\cdot \delta_{(12)}  + 3 \cdot \delta_{(12)(34)} + 8 \cdot \delta_{(123)} + 6 \cdot \delta_{(1234)}}{24} \\
&= \frac{1 + 6\cdot \delta_{(12)} + 9 \cdot \delta_{(12)(34)}  + 8 \cdot \delta_{(123)} }{24}.
\end{align*}
Iterating over all $\delta_{(12)}, \delta_{(12)(34)}, \delta_{(123)} \in \set{0,1}$, we obtain the desired result.
\end{proof}

\begin{lemma} \label{S-A5} If $G \subseteq \GL_2(\ell)$ is a subgroup of type $\Afive$, then $\mathcal{F}(G) \in \set{\frac{1}{60}, \frac{4}{15}, \frac{7}{20}, \frac{5}{12}, \frac{3}{5}, \frac{2}{3}, \frac{3}{4}, 1}$.
\end{lemma}
\begin{proof} Let $\phi: \overline{G} \to A_5$ be an isomorphism and define $\delta_h$ for each $h \in A_5$ as above. The conjugacy classes of $A_5$ are $[()], [(12)(34)], [(123)], [(12345)]$, and $[(12354)]$, of sizes 1, 15, 20, 12, and 12, respectively. Observe that $(12345)^2 = (13524) \in [(12354)]$. Hence, by Lemma \ref{I-lemm}(2,3), we have that $\delta_{(12345)} = \delta_{(12354)}$. Thus by (\ref{S-ex-sum}),
\begin{align*}
\mathcal{F}(G) 
&= \frac{1 \cdot \delta_{()} + 15 \cdot \delta_{(12)(34)} + 20 \cdot \delta_{(123)} + 12 \cdot \delta_{(12345)} + 12 \cdot \delta_{(12354)}}{60} \\
&= \frac{1 + 15 \cdot \delta_{(12)(34)} + 20 \cdot \delta_{(123)} + 24 \cdot \delta_{(12345)}}{60}.
\end{align*}Iterating over all $\delta_{(12)(34)}, \delta_{(123)},\delta_{(12345)} \in \set{0,1}$, we obtain the desired result.
\end{proof}

\section{Proof of Proposition \ref{K-prop}} \label{K-casework}

In this section, we prove the lemmas that are referenced in our proof of Proposition \ref{K-prop}. We start with some observations that will be occasionally useful. As before, $\ell$ denotes an odd prime throughout.

\begin{lemma} \label{I1-lemm}
Each of the following statements holds:
\begin{enumerate}
\item \label{I1-lemm-1} For $\gamma \in \GL_2(\ell)$, we have that $\gamma \in \mathcal{I}_1(\ell)$ if and only if $\det \gamma + 1 = \tr \gamma$.
\item For subgroups $G_1, G_2 \subseteq \GL_2(\ell)$, if $G_1$ is conjugate to $G_2$ in $\GL_2(\ell)$, then $\mathcal{F}_1(G_1) = \mathcal{F}_1(G_2)$.
\end{enumerate}
\end{lemma}
\begin{proof}
\begin{enumerate}
  \item Note that $\gamma \in \mathcal{I}_1(\ell)$ if and only if $1$ is a root of the characteristic polynomial of $\gamma$, which is given by $p_\gamma(x) = x^2 - \tr \gamma \cdot x + \det \gamma$. As $p_\gamma(1) = 1 - \tr \gamma + \det \gamma$, we see that $1$ is a root of the characteristic polynomial of $\gamma$ if and only if $\det \gamma + 1 = \tr \gamma$.
  \item This follows from the fact that a matrix's eigenvalues are invariant under conjugation in $\GL_2(\ell)$. \qedhere
\end{enumerate}
\end{proof}

We now give a lemma that places a restrictive upper bound on $\mathcal{F}_1(G)$, provided that $\abs{G \cap \mathcal{Z}(\ell)}$ is large. Although the bound holds for arbitrary subgroups of $\GL_2(\ell)$, we shall only employ it in the exceptional cases.

\begin{lemma} \label{proj-lem} If $G \subseteq \GL_2(\ell)$ is a subgroup, then 
\[ \mathcal{F}_1(G) \leq \frac{2}{\abs{G \cap \mathcal{Z}(\ell)}} - \frac{1}{\abs{G}}. \]
In particular, if $\abs{G \cap \mathcal{Z}(\ell)} \geq 3$, then $\mathcal{F}_1(G) \leq \frac{2}{3}$.
\end{lemma}
\begin{proof} For a matrix $\gamma \in \GL_2(\ell)$ and scalar $\lambda \in \F_\ell^\times$, the eigenvalues of the product $\lambda  \gamma$ are the eigenvalues of $\gamma$, multiplied by $\lambda$. Thus, since matrices in $\GL_2(\ell)$ have at most two eigenvalues, in each of the $\abs{\overline{G}}$ fibers of the projection $G \to \overline{G}$, there exist at most two matrices contained in $\mathcal{I}_1(\ell)$. In fact, the kernel of $G \to \overline{G}$ contains only a single matrix in $\mathcal{I}_1(\ell)$, the identity matrix. Applying these observations, we obtain
\[ \mathcal{F}_1(G) = \frac{\abs{G \cap \mathcal{I}_1(\ell)}}{\abs{G}} \leq \frac{2(\abs{\overline{G}} - 1) + 1}{\abs{G}} = \frac{2}{\abs{G}/\abs{\overline{G}}} - \frac{1}{\abs{G}}. \]
Finally, notice that $\abs{G \cap \mathcal{Z}(\ell)} = \abs{G}/\abs{\overline{G}}$, by the first isomorphism theorem applied to $\pi: G \twoheadrightarrow \overline{G}$.
\end{proof}

Subgroups $G \subseteq \GL_2(\ell)$ for which $\abs{G \cap \mathcal{Z}(\ell)} \leq 2$ are problematic from the point of view of the previous lemma. In the case of $\abs{G \cap \mathcal{Z}(\ell)} = 2$, we may say something more that will be useful when considering exceptional subgroups. In order to do so, we introduce the following subset of $\overline{G}$,
\begin{equation*} \label{Gbar2}
\overline{G}_2 \coloneqq \set{\overline{\gamma} \in \overline{G} : \text{ the order of } \overline{\gamma} \text{ is two}}.
\end{equation*}

\begin{lemma} \label{proj-lem2} If $G \subseteq \GL_2(\ell)$ is a subgroup for which $\abs{G \cap \mathcal{Z}(\ell)} = 2$, then
\[ \mathcal{F}_1(G) \leq \frac{2\abs{\overline{G}_2} + \abs{\overline{G} \setminus \overline{G}_2}}{\abs{G}}. \]
\end{lemma}
\begin{proof} The scalar group $\mathcal{Z}(\ell)$ is cyclic of order $\ell-1$. Its one (and only) subgroup of order $2$ is $\set{\pm I}$. Hence, since $\abs{G \cap \mathcal{Z}(\ell)} = 2$, we have that $G \cap \mathcal{Z}(\ell) = \set{\pm I}$. Fix $\overline{\gamma} \in \overline{G}$ and note that its fiber under $G \to \overline{G}$ is $\set{\pm \gamma}$. If $\set{\pm \gamma} \subseteq \mathcal{I}_1(\ell)$, then the eigenvalues of $\gamma$ are $1$ and $-1$, so $\gamma$ is conjugate to $\begin{psmallmatrix} 1 & 0 \\ 0 & -1 \end{psmallmatrix}$ in $\GL_2(\ell)$. In particular, the order of both $\gamma$ in $G$ and $\overline{\gamma}$ in $\overline{G}$ is two. We conclude that the fiber of each matrix in $\overline{G}\setminus\overline{G}_2$ contains at most one matrix in $\mathcal{I}_1(\ell)$. Thus, we have the inequality that is claimed in the statement of the lemma.
\end{proof}

The remainder of this section is devoted to proving the lemmas referenced in our proof of Proposition \ref{K-prop}. We proceed case-by-case, considering subgroups of type $\Cs, \Ns, \B, \Afour, \Afive,$ and $\Sfour$ in Dickson's classification. By Lemma \ref{I1-lemm}, we may (and do) make the assumptions described in the paragraph immediately preceding Section \ref{S-CsB-P}.

\subsection{Split Cartan and Borel subgroups} For a subgroup $G \subseteq \mathcal{B}(\ell)$, we define for $i=1$ and $i=4$ the homomorphism $\psi_i : G \to \F_\ell^\times$ given by $\begin{psmallmatrix} a_1 & a_2 \\ 0 & a_4 \end{psmallmatrix} \mapsto a_i$. We observe that $G \cap \mathcal{I}_1(\ell) = \ker \psi_1 \cup \ker \psi_4$ and thus
\begin{equation} \label{K-CsB-ker}
\abs{G \cap \mathcal{I}_1(\ell)} = \abs{\ker \psi_1} + \abs{\ker \psi_4} - \abs{\ker \psi_1 \cap \ker \psi_4}.
\end{equation}

\begin{lemma} \label{K-Cs}  If $G$ is of type $\Cs$ and $\mathcal{F}_1(G) \neq 1$, then $\mathcal{F}_1(G) \leq \frac{1}{2} + \frac{1}{\abs{G}}$.
\end{lemma}
\begin{proof} Let each $\psi_i: G \to \F_\ell^\times$ be defined for $G$ as above. Since $\ker \psi_1 \cap \ker \psi_4 = \set{I}$, we have by (\ref{K-CsB-ker}) that
\begin{equation} \label{K-Cs-1} \abs{G \cap \mathcal{I}_1(\ell)} = \abs{\ker \psi_1} + \abs{\ker \psi_4} - 1. 
\end{equation}
The subgroup of $G$ generated by $\ker \psi_1 \cup \ker \psi_4$ has order $\abs{\ker \psi_1} \abs{\ker \psi_4} $, so  $\abs{\ker \psi_1} \abs{\ker \psi_4} \leq \abs{G}$. Thus,
\begin{equation} \label{K-Cs-2}
 \abs{\ker \psi_1} + \abs{\ker \psi_4} \leq \abs{\ker \psi_1} + \frac{\abs{G}}{\abs{\ker \psi_1}}
 \end{equation}
Now if $\abs{\ker \psi_1} = \abs{\ker \psi_4} = 1$, then $\mathcal{F}_1(G) = \frac{1}{\abs{G}}$ by (\ref{K-Cs-1}) and we are done. So assume, without loss of generality, that $\abs{\ker \psi_1} > 1$. Because $\ker \psi_1 \subseteq G \cap \mathcal{I}_1(\ell)$ and $\mathcal{F}_1(G) \neq 1$, we further have that $\abs{\ker \psi_1} < \abs{G}$. As $\abs{\ker \psi_1}$ is an integer that divides $\abs{G}$ and satisfies the inequalities $1 < \abs{\ker \psi_1} < \abs{G}$, we have that
\begin{equation} \label{K-Cs-3}
\abs{\ker \psi_1} + \frac{\abs{G}}{\abs{\ker \psi_1}} \leq \frac{1}{2}\abs{G}+ 2.
\end{equation}
Now by combining (\ref{K-Cs-1}), (\ref{K-Cs-2}), and (\ref{K-Cs-3}), we conclude that
\[
\mathcal{F}_1(G) = \frac{\abs{G \cap \mathcal{I}_1(\ell)}}{\abs{G}} \leq \frac{\p{\frac{1}{2}\abs{G}+2}-1}{\abs{G}} = \frac{1}{2} + \frac{1}{\abs{G}},
\]
completing the proof.
\end{proof}

\begin{lemma} \label{K-B} If $G$ is of type $\B$ and $\mathcal{F}_1(G) \neq 1$, then $\mathcal{F}_1(G) \leq \frac{1}{2} + \frac{\ell}{\abs{G}}$.
\end{lemma}
\begin{proof} As $G \subseteq \mathcal{B}(\ell)$ and $\ell$ divides $\abs{G}$, we have that $\begin{psmallmatrix} 1 & 1 \\ 0 & 1 \end{psmallmatrix} \in G$. Thus,
\[ \ker \psi_1 \cap \ker \psi_4 = \left\langle\begin{pmatrix} 1 & 1 \\ 0 & 1 \end{pmatrix}\right\rangle = \set{\begin{pmatrix} 1 & a \\ 0 & 1 \end{pmatrix} : a \in \F_\ell}. \]
Hence, by (\ref{K-CsB-ker}), we have that
\begin{equation} \label{K-B-1} 
\abs{G \cap \mathcal{I}(\ell)} = \abs{\ker \psi_1} + \abs{\ker \psi_4} - \ell. 
\end{equation}
The subgroup of $G$ generated by $\ker \psi_1 \cup \ker \psi_4$ has order $\frac{1}{\ell}\abs{\ker \psi_1} \abs{\ker \psi_4} $, so  $\frac{1}{\ell}\abs{\ker \psi_1} \abs{\ker \psi_4} \leq \abs{G}$. Thus,
\begin{align} \label{K-B-2}
\abs{\ker \psi_1} + \abs{\ker \psi_4} - \ell &\leq \abs{\ker \psi_1} + \frac{\abs{G}}{\frac{1}{\ell}\abs{\ker \psi_1}} - \ell \\
&= \ell\p{\frac{1}{\ell}\abs{\ker \psi_1} + \frac{\frac{1}{\ell}\abs{G}}{\frac{1}{\ell}\abs{\ker \psi_1}} - 1}.
\end{align}
Now if $\abs{\ker \psi_1} = \abs{\ker \psi_4} = \ell$, then $\mathcal{F}_1(G) = \frac{\ell}{\abs{G}}$ by (\ref{K-B-1}) and we are done. So assume, without loss of generality, that $\abs{\ker \psi_1} > \ell$. Because $\ker \psi_1 \subseteq G \cap \mathcal{I}_1(\ell)$ and $\mathcal{F}_1(G) \neq 1$, we further have that $\abs{\ker \psi_1} < \abs{G}$. As $\frac{1}{\ell}\abs{\ker \psi_1}$ is an integer that divides $\frac{1}{\ell}\abs{G}$ and satisfies the inequalities $1 < \frac{1}{\ell}\abs{\ker \psi_1} < \frac{1}{\ell} \abs{G}$, we have that
\begin{equation} \label{K-B-3}
\frac{1}{\ell}\abs{\ker \psi_1} + \frac{\frac{1}{\ell}\abs{G}}{\frac{1}{\ell}\abs{\ker \psi_1}} \leq \frac{1}{2}\cdot \frac{1}{\ell} \abs{G} + 2.
\end{equation}
Now by combining (\ref{K-B-1}), (\ref{K-B-2}), and (\ref{K-B-3}), we conclude that
\[
\mathcal{F}_1(G) = \frac{\abs{G \cap \mathcal{I}_1(\ell)}}{\abs{G}} \leq \frac{\ell\p{\p{\frac{1}{2\ell}\abs{G}+2}-1}}{\abs{G}} = \frac{1}{2} + \frac{\ell}{\abs{G}}. \qedhere
\]
\end{proof}

The  above lemma leaves open the possibility of a subgroup $G$ of type $\B$ satisfying $\frac{3}{4} < \mathcal{F}_1(G) < 1$ in the single case when $\frac{1}{\ell}\abs{G} = 3$. As we see in the following remark, this case in fact presents no issues.
\begin{remark} \label{K-B-ext}  If in the above lemma, the quantity $\frac{1}{\ell}\abs{G}$ is prime, then $\mathcal{F}_1(G) = \frac{\ell}{\abs{G}}$. Indeed, for $i=1,4$ we have that $\frac{1}{\ell}\abs{\ker \psi_i}$ divides $\frac{1}{\ell}\abs{G}$ and satisfies the inequalities $1 \leq \abs{\ker \psi_i} \leq \frac{1}{\ell} \abs{G}$. Thus $\abs{\ker \psi_i} \in \set{\ell, \abs{G}}$ for each $i=1,4$. As noted in the proof of the lemma, if $\abs{\ker \psi_1} = \abs{\ker \psi_4} = \ell$, then $\mathcal{F}(G) = \frac{\ell}{\abs{G}}$. The case of $\abs{\ker \psi_1} = \abs{\ker \psi_4} = \abs{G}$ cannot occur since then the inequality $\frac{1}{\ell}\abs{\ker \psi_1} \abs{\ker \psi_4} \leq \abs{G}$ is violated. Both of the remaining cases are excluded by the assumptions of the lemma, since in each we have $\mathcal{F}(G) = \frac{\abs{G} + \ell - \ell}{\abs{G}} = 1$.
\end{remark}

\subsection{Normalizer of the split Cartan subgroup}

\begin{lemma} \label{K-Ns} If $G \subseteq \GL_2(\ell)$ is a subgroup of type $\Ns$ for which $\mathcal{F}_1(G) \neq 1$, then $\mathcal{F}_1(G) \leq \frac{1}{2} + \frac{1}{\abs{G}}$.
\end{lemma}
\begin{proof}  Write $G_c \coloneqq G \cap \mathcal{C}_s(\ell)$ and $G_n \coloneqq G \setminus G_c$. If $G_n \cap \mathcal{I}_1(\ell) = \emptyset$, then we are done as by Lemma  \ref{K-Cs},
\begin{align*}
\mathcal{F}_1(G) &= \frac{\abs{G_c \cap \mathcal{I}_1(\ell)} + \abs{G_n \cap \mathcal{I}_1(\ell)}}{\abs{G}} \\
&\leq \frac{\p{\frac{1}{2} \abs{G_c} + 1} + 0}{\abs{G}} \\
&= \frac{\frac{1}{4}\abs{G}+1}{\abs{G}} \\
&= \frac{1}{4} + \frac{1}{\abs{G}}.
\end{align*}
So we shall assume that $G_n \cap \mathcal{I}_1(\ell) \neq \emptyset$. Say $\gamma_0 \in G_n \cap \mathcal{I}_1(\ell)$ and note that $\tr\gamma_0=0$ and $\tr(\gamma \gamma_0) = 0$. Hence Lemma \ref{I1-lemm}(\ref{I1-lemm-1}), gives that
\[
\gamma \gamma_0 \in \mathcal{I}_1(\ell) \iff \det(\gamma \gamma_0) = -1 \iff \det \gamma = 1.
\]
Thus, we have that
\[ G_n \cap \mathcal{I}_1(\ell) = \set{\gamma \gamma_0 : \gamma \in G_c \cap \SL_2(\ell)}. \]
Either $[G_c : G_c \cap \SL_2(\ell)] \geq 2$ or $G_c = G_c \cap \SL_2(\ell)$. In the former case, we are done as then
\begin{align*}
\mathcal{F}_1(G) &= \frac{\abs{G_c \cap \mathcal{I}_1(\ell)} + \abs{G_n \cap \mathcal{I}_1(\ell)}}{\abs{G}} \\
&\leq \frac{\p{\frac{1}{2} \abs{G_c} + 1} + \frac{1}{2}\abs{G_c}}{\abs{G}} \\
&= \frac{\frac{1}{2}\abs{G}+1}{\abs{G}} \\
&= \frac{1}{2} + \frac{1}{\abs{G}}.
\end{align*}
So we consider the case of $G_c = G_c \cap \SL_2(\ell)$. It is clear that $\mathcal{C}_s(\ell) \cap \SL_2(\ell) \cap \mathcal{I}_1(\ell) = \set{I}$, so in particular $G_c \cap \mathcal{I}_1(\ell) = \set{I}$. Hence, in this case, we also have the bound
\[
\mathcal{F}_1(G) = \frac{\abs{G_c \cap \mathcal{I}_1(\ell)} + \abs{G_n \cap \mathcal{I}_1(\ell)}}{\abs{G}}
= \frac{1 + \abs{G_n}}{\abs{G}} = \frac{1 + \frac{1}{2}\abs{G}}{\abs{G}} = \frac{1}{2} + \frac{1}{\abs{G}},
\]
completing the proof.
\end{proof}

\subsection{Exponential subgroups} We finally consider $G \subseteq \GL_2(\ell)$ an exceptional subgroup. We split our consideration into three cases: $\abs{G \cap \mathcal{Z}(\ell)} = 1$, $\abs{G \cap \mathcal{Z}(\ell)} = 2$, and $\abs{G \cap \mathcal{Z}(\ell)} \geq 3$. It is not clear (at least to the author) whether the first of these three cases may occur, so we pose the following question: Does there exist an exceptional subgroup $G \subseteq \GL_2(\ell)$ for which ${G \cap \mathcal{Z}(\ell)} = \set{I}$?

Lacking an affirmative answer, we start by considering the (possibly vacuous) case of $G \cap \mathcal{Z}(\ell) = \set{I}$. We proceed via conjugacy class considerations, as in Lemmas \ref{S-A4}, \ref{S-S4}, and \ref{S-A5}. First, two quick observations. 

\begin{lemma} \label{K-exp-auxlem} For $\gamma \in \GL_2(\ell)$, we have that
\begin{enumerate}
  \item $\gamma \in \mathcal{I}_1(\ell)$ implies $\gamma^2 \in \mathcal{I}_1(\ell)$, and
  \item $\gamma^2 = I$ implies $\gamma \in \mathcal{I}_1(\ell)$ or $\gamma = - I$.
\end{enumerate}
\end{lemma}
\begin{proof}
\begin{enumerate}
  \item This is clear, since if $1$ is an eigenvalue of $\gamma$, then $1^2$ is an eigenvalue of $\gamma^2$.
  \item Here the minimal polynomial of $\gamma$ divides $X^2 - I$. Thus only if the minimal polynomial of $\gamma$ equals $X + I$ may $\gamma \not\in \mathcal{I}_1(\ell)$. But then $\gamma = - I$. \qedhere
\end{enumerate}
\end{proof}

\begin{lemma} \label{K-exp-lem} If $G \subseteq \GL_2(\ell)$ is an exceptional for which $G \cap \mathcal{Z}(\ell) = \set{I}$ and $\mathcal{F}_1(G) \neq 1$, then $\mathcal{F}_1(G) \leq \frac{3}{4}$.
\end{lemma}
\begin{proof} We have three cases to consider as $G$ may be isomorphic to $A_4, S_4,$ or $A_5$. Below, we write $[\gamma]$ to denote the conjugacy class of a matrix $\gamma$ in $G$.

First assume that $G \cong A_4$. The conjugacy classes of $A_4$ have sizes 1,3,4, and 4. Fix $\gamma \in G$ with $\gamma \not\in\mathcal{I}_1(\ell)$. Then $\gamma \neq I$, so its conjugacy class $[\gamma]$ has size at least 3. We have that $[\gamma] \cap \mathcal{I}_1(\ell) = \emptyset$, so  $\mathcal{F}_1(G) \leq \frac{12-3}{12} = \frac{3}{4}$.

Now assume that $G \cong S_4$. The conjugacy classes of $S_4$ have sizes 1,6,3,8, and 6. Fix $\gamma \in G$ with $\gamma \not\in \mathcal{I}_1(\ell)$. Note that the conjugacy class of size 3 consists of elements of order 2. Thus by Lemma \ref{K-exp-auxlem}(2) and our assumption that $G \cap \mathcal{Z}(\ell) = \set{I}$, we have that the size of $[\gamma]$ is at least 6. Hence $\mathcal{F}_1(G) \leq \frac{24 - 6}{24} = \frac{3}{4}$.

Finally assume that $G \cong A_5$. The conjugacy classes of $A_5$ have sizes 1, 15, 20, 12, and 12. Fix $\gamma \in G$ with $\gamma \not\in \mathcal{I}_1(\ell)$. If the size of $[\gamma]$ is 15 or 20, then we are done as then $\mathcal{F}_1(G) \leq \frac{60-15}{60} = \frac{3}{4}$. So we shall assume that $[\gamma]$ is one of the conjugacy classes of size $12$. Let $\gamma_0 \in G$ be such that $[\gamma_0]$ is the other conjugacy class of $G$ of size 12. Then $\gamma_0^2 \in [\gamma]$, so by the contrapositive of Lemma \ref{K-exp-auxlem}(1), we have that $\gamma_0 \not\in \mathcal{I}_1(\ell)$. Thus $\p{[\gamma] \cup [\gamma_0]} \cap \mathcal{I}_1(\ell) = \emptyset$. Consequently, $\mathcal{F}_1(G) \leq \frac{60-2\cdot12}{60} = \frac{3}{4}$.
\end{proof}

Next, we consider the case of $\abs{G \cap \mathcal{Z}(\ell)} = 2$. Here we proceed via Lemma \ref{proj-lem2}.

\begin{lemma} \label{K-exp} If $G \subseteq \GL_2(\ell)$ is an exceptional for which $\abs{G \cap \mathcal{Z}(\ell)} = 2$, then
\[
\mathcal{F}_1(G) \leq 
\begin{cases}
\frac{5}{8} &  G \text{ is of type } \Afour \text{ or } \Afive \\
\frac{11}{16} & G \text{ is of type } \Sfour.
\end{cases}
\]
\end{lemma}
\begin{proof} First assume $G$ is of type $\Afour$. The group $A_4$ has 12 elements, of which 3 have order 2. Thus, by Lemma \ref{proj-lem2}, we have that $\mathcal{F}_1(G) \leq \frac{2 \cdot 3 + (12-3)}{24} = \frac{5}{8}$. We obtain the upper bounds for the other cases similarly. Specifically, apply Lemma \ref{proj-lem2} on noting that $A_5$ has 60 elements, of which 15 have order 2 and that $S_4$ has 24 elements, of which 9 have order 2. 
\end{proof}

Finally, we note that the case of $\abs{G \cap \mathcal{Z}(\ell)} \geq 3$ has already been handled via Lemma \ref{proj-lem}.

\section{Densities for Non-CM Elliptic Curves over the Rationals} \label{appendix}

Let $E/\Q$ be an elliptic curve without complex multiplication. For a prime number $\ell$, let $G_E(\ell)$ denote the image of the mod $\ell$ Galois representation $\rho_{E,\ell}: \Gal(\overline{\Q}/\Q) \to \GL_2(\ell)$. Serre's open image theorem \cite[Th\'{e}or\`{e}m 3]{Se1972} gives that $G_E(\ell) = \GL_2(\ell)$ for all sufficiently large $\ell$. If $G_E(\ell) = \GL_2(\ell)$, then
\[
\delta(\mathcal{S}^1_{E,\ell}) = \frac{\ell^2-2}{(\ell^2-1)(\ell-1)} \quad \text{and} \quad \delta(\mathcal{S}_{E,\ell}) = \frac{\ell+2}{2(\ell+1)},
\]
as we see by Proposition \ref{bridge-prop} and a calculation of $\mathcal{F}_1(\GL_2(\ell))$ and $\mathcal{F}(\GL_2(\ell))$ (the latter is carried out in Lemma \ref{S-PGL} and the former follows similarly). 

A prime $\ell$ is \textit{exceptional} for $E$ if $G_E(\ell) \neq \GL_2(\ell)$ and, in this instance, the group $G_E(\ell)$ is called an \textit{exceptional image} for $\ell$. All exceptional images are known \cite{Zy2015} for $\ell \leq 11$. For primes $\ell \geq 13$, as a result of systematic computations  \cite{Su2016} and significant partial results (e.g. \cite{Ba2019},\cite{Bi2011}, \cite{Bi2013},  \cite{Ma1977}, \cite{Ma1978}, \cite{Se1972}), it is conjectured that all exceptional images are known and that $G_E(\ell) = \GL_2(\ell)$ for $\ell > 37$.

We reproduce from \cite[Table 3]{Su2016} the conjecturally complete list of 63 exceptional images in the first column of the table below. For each exceptional image $G$, we list in columns two and three the associated values of $\delta(\mathcal{S}^1_{E,\ell})$ and $\delta(\mathcal{S}_{E,\ell})$ for elliptic curves $E/\Q$ with $G_E(\ell) = G$. These densities are straightforward and fast to compute as, by Proposition \ref{bridge-prop}, we simply need to compute the proportions $\mathcal{F}_1(G)$ and $\mathcal{F}(G)$.

\renewcommand{\arraystretch}{1.1}
\begin{table}[h] 
\begin{tabular}{lcc} 
$G_E(\ell)$ & $\delta(\mathcal{S}^1_{E,\ell})$ & $\delta(\mathcal{S}_{E,\ell})$  \\ \hline
$\la{2Cs}$ & $1$ & $1$ \\
$\la{2B}$ & $1$ & $1$ \\
$\la{2Cn}$ & $\frac{1}{3}$ & $\frac{1}{3}$ \\
$\la{3Cs.1.1}$ & $1$ & $1$ \\
$\la{3Cs}$ & $\frac{3}{4}$ & $1$ \\
$\la{3B.1.1}$ & $1$ & $1$ \\
$\la{3B.1.2}$ & $1$ & $1$ \\
$\la{3Ns}$ & $\frac{5}{8}$ & $\frac{3}{4}$ \\
$\la{3B}$ & $\frac{3}{4}$ & $1$ \\
$\la{3Nn}$ & $\frac{5}{16}$ & $\frac{3}{8}$ \\
$\la{5Cs.1.1}$ & $1$ & $1$ \\
$\la{5Cs.1.3}$ & $\frac{1}{2}$ & $1$ \\
$\la{5Cs.4.1}$ & $\frac{5}{8}$ & $1$ \\
$\la{5Ns.2.1}$ & $\frac{3}{16}$ & $\frac{1}{2}$ \\
$\la{5Cs}$ & $\frac{7}{16}$ & $1$ \\
$\la{5B.1.1}$ & $1$ & $1$ \\
$\la{5B.1.2}$ & $1$ & $1$ \\
$\la{5B.1.4}$ & $\frac{1}{2}$ & $1$ \\
$\la{5B.1.3}$ & $\frac{1}{2}$ & $1$ \\
$\la{5Ns}$ & $\frac{11}{32}$ & $\frac{3}{4}$ \\
$\la{5B.4.1}$ & $\frac{5}{8}$ & $1$ \\
\end{tabular}
\;
\begin{tabular}{lcc}
$G_E(\ell)$ & $\delta(\mathcal{S}^1_{E,\ell})$ & $\delta(\mathcal{S}_{E,\ell})$  \\ \hline
$\la{5B.4.2}$ & $\frac{5}{8}$ & $1$ \\
$\la{5Nn}$ & $\frac{7}{48}$ & $\frac{1}{3}$ \\
$\la{5B}$ & $\frac{7}{16}$ & $1$ \\
$\la{5S4}$ & $\frac{19}{96}$ & $\frac{5}{12}$ \\
$\la{7Ns.2.1}$ & $\frac{4}{9}$ & $1$ \\
$\la{7Ns.3.1}$ & $\frac{11}{36}$ & $1$ \\
$\la{7B.1.1}$ & $1$ & $1$ \\
$\la{7B.1.3}$ & $1$ & $1$ \\
$\la{7B.1.2}$ & $\frac{1}{3}$ & $1$ \\
$\la{7B.1.5}$ & $\frac{1}{3}$ & $1$ \\
$\la{7B.1.6}$ & $\frac{2}{3}$ & $1$ \\
$\la{7B.1.4}$ & $\frac{2}{3}$ & $1$ \\
$\la{7Ns}$ & $\frac{17}{72}$ & $\frac{3}{4}$ \\ 
$\la{7B.6.1}$ & $\frac{7}{12}$ & $1$ \\
$\la{7B.6.3}$ & $\frac{7}{12}$ & $1$ \\
$\la{7B.6.2}$ & $\frac{1}{4}$ & $1$ \\
$\la{7Nn}$ & $\frac{3}{32}$ & $\frac{5}{16}$ \\
$\la{7B.2.1}$ & $\frac{4}{9}$ & $1$ \\
$\la{7B.2.3}$ & $\frac{4}{9}$ & $1$ \\
$\la{7B}$ & $\frac{11}{36}$ & $1$ \\
$\la{11B.1.4}$ & $\frac{1}{5}$ & $1$ \\
\end{tabular} 
\;
\begin{tabular}{lcc}
$G_E(\ell)$ & $\delta(\mathcal{S}^1_{E,\ell})$ & $\delta(\mathcal{S}_{E,\ell})$  \\ \hline
$\la{11B.1.6}$ & $\frac{1}{5}$ & $1$ \\
$\la{11B.1.5}$ & $\frac{1}{5}$ & $1$ \\
$\la{11B.1.7}$ & $\frac{1}{5}$ & $1$ \\
$\la{11B.10.4}$ & $\frac{3}{20}$ & $1$ \\
$\la{11B.10.5}$ & $\frac{3}{20}$ & $1$ \\
$\la{11Nn}$ & $\frac{13}{240}$ & $\frac{7}{24}$ \\
$\la{13S4}$ & $\frac{35}{288}$ & $\frac{3}{4}$  \\
$\la{13B.3.1}$ & $\frac{7}{18}$ & $1$ \\
$\la{13B.3.2}$ & $\frac{7}{18}$ & $1$ \\
$\la{13B.3.4}$ & $\frac{2}{9}$ & $1$ \\
$\la{13B.3.7}$ & $\frac{2}{9}$ & $1$ \\
$\la{13B.5.1}$ & $\frac{5}{16}$ & $1$ \\
$\la{13B.5.2}$ & $\frac{5}{16}$ & $1$ \\
$\la{13B.5.4}$ & $\frac{7}{48}$ & $1$ \\
$\la{13B.4.1}$ & $\frac{17}{72}$ & $1$ \\
$\la{13B.4.2}$ & $\frac{17}{72}$ & $1$ \\
$\la{13B}$ & $\frac{23}{144}$ & $1$ \\
$\la{17B.4.2}$ & $\frac{11}{64}$ & $1$ \\
$\la{17B.4.6}$ & $\frac{11}{64}$ & $1$ \\
$\la{37B.8.1}$ & $\frac{47}{432}$ & $1$ \\
$\la{37B.8.2}$ & $\frac{47}{432}$ & $1$ \\
\end{tabular} 
\end{table}

\end{document}